\documentclass[12pt]{article}
\usepackage{amssymb}
\usepackage{color}
\usepackage{amsmath}
\usepackage{graphicx}
\usepackage{theorem}
\usepackage{here}

\newtheorem{thm}{Theorem}[section]
\newtheorem{lem}{Lemma}[section]
\newtheorem{lema}{Lemma A.}

\theorembodyfont{\rmfamily}

\newtheorem{remark}{Remark}[section]
\makeatletter
\renewcommand{\theequation}{%
   \thesection.\arabic{equation}}
\@addtoreset{equation}{section}
\makeatother

\def\bde{{\text{\boldmath $\delta$}}}

\def\bmu{{\text{\boldmath $\mu$}}}

\def\bdeh{{\widehat \bde}}

\def\De{{\Delta}}

\def\Deh{{\widehat \De}}

\def\x{{\text{\boldmath $x$}}}
\def\y{{\text{\boldmath $y$}}}
\def\z{{\text{\boldmath $z$}}}

\def\E{{\mathrm E}}

\def\xb{{\overline \x}}

\def\ah{{\hat a}}

\def\ch{{\hat c}}

\def\Uh{{\widehat U}}
\def\Vh{{\widehat V}}

\def\Nc{{\cal N}}

\def\tr{{\rm tr\,}}

\def\[{{\text{\boldmath $[$}}}
\def\]{{\text{\boldmath $]$}}}

\def\0{{\bf 0}}
\def\1{{\bf 1}}

\def\/{{\Bigr/\!\!}}

\def\1r{{\rm (1)}}
\def\2r{{\rm (2)}}
\def\3r{{\rm (3)}}
\def\4r{{\rm (4)}}
\def\5r{{\rm (5)}}

\begin{document}
\title{Asymptotic Properties of the Misclassification Errors for Euclidean Distance Discriminant Rule in High-Dimensional Data}
\author{
Hiroki Watanabe
\footnote{Department of Mathematical Information Science}, 
{Masashi Hyodo\footnote{Department of Mathematical Information Science}, 
Takashi Seo\footnote{Department of Mathematical Information Science} 
and 
Tatjana Pavlenko
\footnote{Department of Mathematics, KTH Royal Institute of Technology}
}
\\
{\it\normalsize Tokyo University of Science and KTH Royal Institute of Technology}
}
\maketitle

\begin{abstract}
Performance accuracy of the Euclidean Distance Discriminant rule (EDDR) is
studied in the high-dimensional asymptotic framework which allows the dimensionality to
exceed sample size. Under mild assumptions on the traces of the
covariance matrix, our new results provide the asymptotic distribution of the conditional
misclassification error and the explicit expression for the consistent and
asymptotically unbiased estimator of the expected misclassification error. To get
these properties, new results on the asymptotic normality of the quadratic forms
and traces of the higher power of Wishart matrix, are established.  Using  our asymptotic results,
we further develop two  generic methods of determining a cut-off point for EDDR  to adjust the
misclassification errors.  Finally, we numerically justify the high accuracy of our asymptotic
findings along with the cut-off determination methods in finite sample applications, inclusive of the large sample
and high-dimensional scenarios.

\par\vspace{4mm}
\end{abstract}

\section{Introduction}
\
\medskip
In this paper, we focus on the discrimination problem which is concerned with the 
allocation of a given object, $\x$, a random vector represented by a set of features $(x_{1},\dots,x_{p})$, to one ot two populations, $\Pi_1$ and $\Pi_2$ given by $\Nc_p(\bmu_1, \Sigma)$ and $\Nc_p(\bmu_2, \Sigma)$, respectively,  where 
$\bmu_1\neq \bmu_2$ and common covariance matrix $\Sigma$ is non-singular.  Let $\{ \x_{gj}\}^{N_g}_{j=1}$ be a random sample of independent observations drawn from $g$th population $\Nc_p(\bmu_g, \Sigma)$, $g=1, 2$.
Let also $N=N_{1}+N_{2}$ denote the total sample size and set $n=N-2$. We are interested to explore the discrimination procedure that can accomodate $p>n$ cases, with the main focus on the performance accuracy in the asymptototic framework that allows $p$ to grow together with $n$.


Clearly, the classical discriminant procedures, like Fisher linear discriminant rule, cannot be used when $p>n$ since the sample covariance matrix is singular and hence cannot be inverted. An intuitively appealing alternative considered in this study focuses on geometrical properties of the sample space and re-formulates the classification problem in terms of the {\it Euclidian distance discriminant rule} (EDDR): assign  a new observation $\x$ to the "nearest" population $\Pi_{g}$, i.e. assign to $\Pi_{g}$ if it is on average closer to the data from $\Pi_{g}$ than to the data from the other population.  Matusita's papers 
(see Matusita (1955), and Matusita and Motoo (1956)) are perhaps the oldest references dealing  with the discriminant rule based on distance measures,  including the case when the multivariate distributions underlying the data are not specified. 

Recently, Aoshima and Yata (2013) have been considered the EDDR for the high-dimensional multi-class problem with 
different class covariance matrices. In particular, they derived asymptotic conditions which ensure that the expected misclassification error converges to zero. Recent paper by Srivastava (2006) used the Moore-Penrose inverse of the estimated covariance matrix and suggested a 
second-order approximation of the expected error rate in high-dimensional data. 

We, in this study, focus on the asymptotic behavior of the misclassification errors of EDDR.
Continuing with the normality assumption, with $\bmu_g$ acting as the centered of the $\Pi_{g}$'s distribution 
we define 
\begin{eqnarray}
T_0(\boldsymbol{x})=\parallel\boldsymbol{x}-\bmu_2\parallel^2-
\parallel\boldsymbol{x}-\bmu_1\parallel^2, 
\label{eqn:t_0}
\end{eqnarray}
and its sample based version as  
\begin{eqnarray}
\widetilde{T}(\x) = \parallel\boldsymbol{x}-\overline{\boldsymbol{x}}_2\parallel^2-
\parallel\boldsymbol{x}-\overline{\boldsymbol{x}}_1\parallel^2 
\label{eqn:ttil}
\end{eqnarray}
where $\| \cdot \|$ denotes the Euclidean norm and $\overline{\x}_g$'s denote the sample mean vectors, $g=1,2$. Hence, each term in 
$(1.1)$ and $(1.2)$ represents the distance between the observed vector $\x$ and the centroid of 
$\Pi_{g}$'s or its sample based counterpart.

The natural advantage of using $\widetilde{T}(\x)$ for classifying high-dimensional data is its ability to mitigate the effect of dimensionality on the performance accuracy. 
Indeed, as it is seen from $(\ref{eqn:ttil})$, $\widetilde{T}(\x)$ utilizes only the marginal distribution of the $p$ variables, thereby 
naturally reducing the effect of large $p$ in implementations. 
But the dimensionality has impact on the classification accuracy. To show this, we first point out that classifier 
$\widetilde{T}(\x)$ has a bias. In fact, 
\begin{eqnarray*}
\E[\widetilde{T}(\x) | \x \in \Pi_g] &=&(-1)^{g-1}\|\bmu_1-\bmu_2\|^2
+\frac{N_1 - N_2}{N_1 N_2} \tr{\Sigma},~g=1,2, 
\end{eqnarray*}
and thus the impact of dimensionality is implied by the quantity $(N_1 - N_2) \tr\Sigma/(N_1N_2)$.  
In this study,  we introduce the bias-corrected version $\widetilde{T}(\x)$ defined as 
\begin{eqnarray}
T(\x) = \|\x - \overline{\x}_2\|^2 -
\|\x - \overline{\x}_1\|^2
- \frac{N_1-N_2}{N_1 N_2} \tr S,
\label{eqn:tt}
\end{eqnarray}
where the subtraction of $(N_1-N_2)/(N_1 N_2)\tr S$ in $(\ref{eqn:tt})$ is to guarantee that $\E[T(\x)| \x\in \Pi_g]=(-1)^{g-1}\|\bmu_1-\bmu_2 \|^2,~g=1,2.$ Here, $S= (1/n) \sum_{g=1}^2 \sum_{j=1}^{N_g}(\x_{gj}-\overline{\x}_g)(\x_{gj}-\overline{\x}_g)^\prime.$

Now, the EDDR given by $T(\x)$ places a new observation $\x$ to $\Pi_1$ if $T(\x) > \tilde{c}$, and to $\Pi_2$ otherwise, 
where $\tilde{c}$ is an appropriate cut-off point. Then, for a specific $\tilde{c}$, the performance accuracy 
of EDDR will be represented by the pair of misclassification error rates that result. 
Precisely, we define the conditional misclassification error of EDDR by 
\begin{eqnarray*}
ce(2|1) =\Pr(T(\x) \leq \tilde{c} |\x\in \Pi_1,~\overline{\x}_1,~\overline{\x}_2,~S)
\end{eqnarray*}
and its expected version by $e(2|1)=\E[ce(2|1)]$, where the expectation is taken with respect to $\overline{\x}_1$, $\overline{\x}_2$ and $S$. 
Our main objective is to derive characteristic properties of both conditional and expected misclassification error
in high-dimensional data.

In many practical problems one type of misclassification error 
is generally regarded as more serious than the other, examples include e.g. medical applications 
associated with the diagnosis of diseases. In such a  case, it might be desired to determine the cut-off $\tilde{c}$ to obtain a specified probability of the error, or at least to approximate a specified probability.  Then, one might base the choice of $\tilde{c}$ on the expected misclassification error.  This method, denoted in what follows  by $\mathbf{M1}$, suggests to set a cut-off point $\tilde{c}$ 
such that 
\begin{eqnarray*}
\mathbf{M1}:  \quad   e(2|1)=\E[ce(2|1)] = \alpha,
\label{M1}
\end{eqnarray*}
where $\alpha$ is a value given by experimenters.


On the other hand,  one may exploit the confidence of the conditional error rate when determining $\tilde{c};$ we denote this method by
\begin{eqnarray*}
\mathbf{M2}:  \quad\Pr(ce(2|1) <eu)=1-\beta,
\label{M2}
\end{eqnarray*}
where $1-\beta$ is the desired level of confidence and $eu$ is an upper bound. 

Both determination methods $\mathbf{M1}$ and  $\mathbf{M2}$ have been established by using large sample approximation, see Anderson (1973), 
McLachlan (1977) and Shutoh et al. (2012). In this study, we extend the consideration to the high-dimensional case. 
Our main theoretical results provide the asymptotically unbiased and consistent estimator of $e(2|1)$ and the limit distribution of $ce(2|1)$ under general assumptions covering the case when $p > n$. In fact, ${\bf M1}$ and ${\bf M2}$ 
procedures can be considered as specific examples of using our generic results in the theory of EDDR in high-dimensions. 

The remaining part of the paper is organized as follows. 
In Section 2, we derived the asymptotically unbiased and consistent estimator of $e(2|1)$. 
Further, the limiting approximations of the cut-off point defined by $\mathbf{M1}$ are established by using this estimator. In Section 3, two estimators of the confidence-based cut-off point
defined by $\mathbf{M2}$ are proposed, for which the asymptotic normality of the conditional error rate is shown. 
Section 4 summaries the results of numerical experiments justifying the validity of the suggested cut-off
estimators for various strength of dependence underlying the data along with a number of high-dimensional scenarios 
where $p$ far exceeds the sample size. We conclude in Section 5, and give a through proofs of newly established 
asymptotic results together with some auxiliary lemmas in Appendix A.     

\section{Evaluation of the expected misclassification error}          %

Getting the closed-form expression for the 
expected error is too demanding, therefore we first shall derive its asymptotic approximation, and then based on this result, 
propose the consistent and asymptotically unbiased estimator of $e(2|1)$ in high dimensions. 
We further show how these results can be used to provide the cut-off by the determination procedure $\mathbf{M1}$. 


Let $\bde=\bmu_1-\bmu_2, a_i=\tr\Sigma^i/p,i=1,\dots,8,\Delta_i=\bde'\Sigma^i\bde,i=1,\dots,7$ and $\Delta_0=\bde'\bde.$
We make the following assumptions for the consistency and unbiasedness of the estimator of $e(2|1)$:

\medskip
$
\displaystyle ({\rm A1}):0<\lim_{(n,p) \to \infty}\frac{p}{n}<\infty, ~~~~0<\lim_{(n,N_i) \to \infty}\frac{N_i}{n+2}<1,~i=1,2.
$

\medskip

$
\displaystyle ({\rm A2}):0<\lim_{(n,p) \to \infty}\Delta_0, \lim_{(n,p) \to \infty}\Delta_1<\infty, ~~~~
0<\lim_{(n,p) \to \infty}a_1, \lim_{(n,p) \to \infty}a_2<\infty.
$

\medskip

$
\displaystyle ({\rm A3}):\lim_{(n,p) \to \infty}\frac{\Delta_3}{n} \to 0, \lim_{(n,p) \to \infty}\frac{a_4}{n} \to 0.
$

\noindent
Assume henceforth $\x\in\Pi_1$. The symmetry of our classification rule makes the probability of error if the mean of $\x$ is
$\bmu_1$ the same as that under $\bmu_2$.  Then for the conditional distribution of $T(\x)$ given $(\xb_1,\xb_2,S)$ it holds that 
\begin{eqnarray*}
T(\x)| (\xb_1,\xb_2,S) \sim \Nc\left( -2U-\frac{N_1-N_2}{N_1N_2}\tr S, 4V\right),
\end{eqnarray*}
where
\begin{align*}
U=& (\xb_1-\xb_2)'(\xb_1-\bmu_1)-\frac{1}{2}(\xb_1-\xb_2)'(\xb_1-\xb_2),\\
V=& (\xb_1-\xb_2)'\Sigma (\xb_1-\xb_2).
\end{align*}
Now the expected error rate $e(2|1)$ of $T(\x)$ can be expressed in terms of $U$ and $V$ as  
\begin{eqnarray}
e(2|1)=\E[ce(2|1)]=\E\left[\Phi\left(\frac{U+(N_2^{-1}-N_1^{-1})p\hat{a}_1/2+c}{\sqrt{V}}\right)\right],
\end{eqnarray}
where the expectation is with respect to $U$ and $V$, $c=\tilde{c}/2,$ $\ah_1=\tr S/p$ and $\Phi(\cdot)$ is the cumulative distribution function of 
the standard normal distribution.  

In order to proceed to asymptotic approximation of $e(2|1)$, 
we need some preparatory stochastic evaluation of $U$ and $V$. 
We introduce the auxiliary random variables  
\begin{eqnarray*} 
\textrm{\boldmath $z$}_1 &=&
N^{- \frac{1}{2}}\Gamma^\prime\Sigma^{- \frac{1}{2}}
(N_1\overline{\textrm{\boldmath $x$}}_1 
+ N_2 \overline{\textrm{\boldmath $x$}}_2 
- N_1 \textrm{\boldmath $\mu$}_1 
- N_2 \textrm{\boldmath $\mu$}_2),\\
\textrm{\boldmath $z$}_2 &=&
\left( \frac{N}{N_1 N_2} \right)^{- \frac{1}{2}}\Gamma^\prime\Sigma^{- \frac{1}{2}}
(\overline{\textrm{\boldmath $x$}}_1
- \overline{\textrm{\boldmath $x$}}_2
- \textrm{\boldmath $\mu$}_1
+ \textrm{\boldmath $\mu$}_2), 
\end{eqnarray*}
and observe that $\textrm{\boldmath $z$}_1$ and $\textrm{\boldmath $z$}_2$ are independent and identically distributed as $\mathcal{N}_p(\textrm{\boldmath $0$},I_p)$, 
where $\Gamma$ is an orthogonal matrix such that 
$\Sigma=\Gamma\Lambda\Gamma'$ and $\Lambda$ is a diagonal matrix of eigenvalues of $\Sigma$. 
By means of $\z_1$ and $\z_2$, we further define
\begin{eqnarray}
U_0&=&-\Delta_0/2, \\
U_1&=&\frac{1}{\sqrt{N}}\textrm{\boldmath $\delta$}^\prime\Gamma\Lambda^{\frac{1}{2}}\textrm{\boldmath $z$}_1
-\left(\frac{N_1}{N N_2} \right)^{\frac{1}{2}}\textrm{\boldmath $\delta$}^\prime\Gamma\Lambda^{\frac{1}{2}}\textrm{\boldmath $z$}_2
+\frac{1}{(N_1 N_2)^{\frac{1}{2}}}\textrm{\boldmath $z$}_1^\prime\Lambda\textrm{\boldmath $z$}_2
-\frac{N_1 -N_2}{2 N_1 N_2}(\textrm{\boldmath $z$}_2^\prime\Lambda\textrm{\boldmath $z$}_2-pa_1), \nonumber\\
& &\\
U_2&=&\frac{(N_1-N_2)p}{2N_1N_2}(\hat{a}_1-a_1), 
\end{eqnarray}
and observe that by using (2.2)-(2.4) the numerator in (2.1) can be decomposed as
\begin{eqnarray}
U+\frac{(N_1-N_2)p\hat{a}_1}{2N_1N_2}=U_0+U_1+U_2.
\end{eqnarray}
By analogy with $U$, $V$ can also be decomposed by first defining $V_0$ and $V_1$ as 
\begin{eqnarray}
V_0 &=&\Delta_1+\frac{Npa_2}{N_1N_2},
~V_1=2 \left(\frac{N}{N_1 N_2} \right)^{\frac{1}{2}} \textrm{\boldmath $\delta$}^\prime \Gamma\Lambda^{\frac{3}{2}}\textrm{\boldmath $z$}_2
+ \frac{N}{N_1 N_2}(\textrm{\boldmath $z$}_2^\prime\Lambda^2\textrm{\boldmath $z$}_2-pa_2)
\end{eqnarray}
and then observing that $V=V_0+V_1$. 
Now for the first moments, we have by (2.5) and (2.6)
\begin{eqnarray*}
\E\left[U + \frac{(N_1-N_2)p\ah_1}{2N_1N_2}\right]=U_0,~
\E[V]=V_0.
\end{eqnarray*}
To evaluate the second moments, we apply Lemma A.3 (see Appendix) and obtain
\begin{eqnarray*}
\E\left[ \left(U+\frac{(N_1-N_2)p\ah_1}{2N_1N_2}-U_0 \right)^2\right] =H_U(\Delta_1, a_2)+o(n^{-1}),
~\E[(V-V_0)^2] =H_V(\Delta_3, a_4),
\end{eqnarray*}
where
\begin{align*}
H_U(\De_1,a_2)=\frac{1}{N_2}\Delta_1
+\frac{(N_1^2 + N_2^2) pa_2}{2 N_1^2 N_2^2},~
H_V(\De_3,a_4)=\frac{4N}{N_1 N_2} \Delta_3
+ \frac{2 N^2pa_4}{(N_1 N_2)^2}.
\end{align*}
Under the assumptions (A1)-(A3), it holds that
\begin{eqnarray}
\E\left[ \left(U+\frac{(N_1-N_2)p\ah_1}{2N_1N_2}-U_0 \right)^2\right] \rightarrow 0,~
\E[(V-V_0)^2] \rightarrow 0,
\end{eqnarray}
and by Chebyshev's inequality, (2.7) implies that
\begin{eqnarray}
U+\frac{(N_1-N_2)p\ah_1}{2N_1N_2} \xrightarrow{P} U_0,~~ V \xrightarrow{P} V_0,
\end{eqnarray}
where $\xrightarrow{P}$ denotes convergence in probability.

Since $\Phi(\cdot)$ in (2.1) is a continuous function of $U$ and $V$, it follows from (2.8), 
by the continuous mapping theorem, that
\begin{eqnarray}
&&\left|\Phi\left(\frac{U+(N_1-N_2)p\ah_1/(2N_1N_2)+c}{\sqrt{V}}\right) - \Phi\left(\frac{U_0+c}{\sqrt{V_0}}\right)\right| \xrightarrow{P} 0. \nonumber
\end{eqnarray}
On the other hand, it naturally holds that
\begin{eqnarray*}
\left|\Phi\left(\frac{U+(N_1-N_2)p\ah_1/(2N_1N_2)+c}{\sqrt{V}}\right) - \Phi\left(\frac{U_0+c}{\sqrt{V_0}}\right)\right| < 1.
\end{eqnarray*}
Hence, by the dominated convergence theorem we have 
\begin{eqnarray}
\E\left[\left|\Phi\left(\frac{U+(N_1-N_2)p\ah_1/(2N_1N_2)+c}{\sqrt{V}}\right) - \Phi\left(\frac{U_0+c}{\sqrt{V_0}}\right)\right|\right] \to 0.
\end{eqnarray}
Further, by applying the Jensen's inequality to (2.9) we get
\begin{eqnarray*}
&&\left|\E\left[\Phi\left(\frac{U+(N_1-N_2)p\ah_1/(2N_1N_2)+c}{\sqrt{V}}\right)\right] - \Phi\left(\frac{U_0+c}{\sqrt{V_0}}\right)\right|\\
&&\leq\E\left[\left|\Phi\left(\frac{U+(N_1-N_2)p\ah_1/(2N_1N_2)+c}{\sqrt{V}}\right) - \Phi\left(\frac{U_0+c}{\sqrt{V_0}}\right)\right|\right]\to 0.
\end{eqnarray*}
The above results  are summarized in the following lemma.
\begin{lem}
Under assumptions {\rm (A1)-(A3)}
\begin{eqnarray}
e(2|1) \to \Phi\left(\frac{U_0+c}{\sqrt{V_0}}\right),
\end{eqnarray}
where $U_0$ and $V_0$ are defined in {\rm (2.2)} and {\rm (2.6)}, respectively.
\end{lem}
In words, Lemma 2.1 provides a closed form expression for the limiting term of $e(2|1)$. 
Hence, to identify the cut-off point for $T(\x)$, we derive a consistent and unbiased estimator of $e(2|1)$ 
by plugging-in consistent estimators of $U_0$ and $V_0$ into the right hand side of (2.10). 

As $U_0$ and $V_0$ are functions of $\De_0$, $\De_1$ and $a_2$, we begin by obtaining their consistent estimators. 
\begin{lem}
Let estimators of $\Delta_0, \Delta_1, a_2$ be defined as 
\begin{eqnarray}
\Deh_0&=&\bdeh^\prime\bdeh-\frac{Np}{N_1N_2}\ah_1,\\
\Deh_1&=&\bdeh^\prime S \bdeh-\frac{Np}{N_1N_2}\ah_2,\\
\hat{a}_2&=&\frac{n^2}{p(n+2)(n-1)} \left ( \tr{S^2} -\frac{(\tr{S})^2}{n} \right ),
\end{eqnarray}
respectively, where $\bdeh=\overline{\textrm{\boldmath $x$}}_1 - \overline{\textrm{\boldmath $x$}}_2.$ Then under assumptions {\rm (A1)-(A3)} 
\begin{eqnarray*}
\Deh_0 \xrightarrow{P} \De_0, ~~\Deh_1 \xrightarrow{P} \De_1,~~ \ah_2 \xrightarrow{P} a_2.
\end{eqnarray*}
\end{lem}
\noindent
{\bf (Proof)}

To show consistency of $a_1$ and $a_2$, we use exact expressions for the variances of these estimators derived in Srivastava (2005) as
\begin{eqnarray}
\E[(\ah_1 - a_1)^2] &=& \frac{2a_2}{np}, \\
\E[(\ah_2 - a_2)^2] &=& \frac{8(n+2)(n+3)(n-1)^2}{p n^5} a_4 +\frac{4(n+2)(n-1)}{n^4}(a_2^2 - p^{-1}a_4).
\end{eqnarray}
Then by applying Chebyshev's inequality to (2.14) and (2.15) it can be seen that
\begin{eqnarray}
\ah_1 \xrightarrow{P} a_1, \ah_2 \xrightarrow{P} a_2.
\end{eqnarray}

To show consistency of $\Deh_0$ and $\Deh_1$, 
we first consider the following random variables
\begin{eqnarray*}
\widetilde{\Delta}_0=\bdeh^\prime\bdeh-\frac{Np}{N_1N_2}a_1,~\widetilde{\Delta}_1=\bdeh^\prime S \bdeh-\frac{Np}{N_1N_2}a_2
\end{eqnarray*}
and evaluate the first two moments of $\bdeh^\prime\bdeh$ and $\bdeh^\prime S \bdeh$. We rewrite
\begin{eqnarray*}
\bdeh^\prime \bdeh&=&
\bde'\bde+2\left(\frac{N}{N_1N_2}\right)^{1/2}\bde'\Sigma^{1/2}\z
+\frac{N}{N_1N_2}\z'\Sigma\z,
\end{eqnarray*}
and
\begin{eqnarray}
\bdeh^\prime S\bdeh&=&
\bde'S\bde+2\left(\frac{N}{N_1N_2}\right)^{1/2}\bde'S\Sigma^{1/2}\z
+\frac{N}{N_1N_2}\z'\Sigma^{1/2}S\Sigma^{1/2}\z,
\end{eqnarray}
where $\z\sim\Nc(\boldsymbol{0},I_p)$. 
Then it easily follows that
\begin{eqnarray}
\E[\widetilde{\Delta}_0]=\Delta_0,~\E[\widetilde{\Delta}_1]=\Delta_1
\end{eqnarray}
and
\begin{eqnarray}
{\rm Var}[\widetilde{\Delta}_0]&=&
\frac{4N}{N_1N_2}\De_1+\frac{2N^2p}{N_1^2N_2^2}a_2, \\
{\rm Var}[\widetilde{\Delta}_1]&=&\frac{2 a_2^2 N^2 p^2}{n N_1^2 N_2^2}+\frac{4 a_2 \De_1 N p}{n N_1 N_2}+\frac{2 a_4 N^3
   p}{n N_1^2 N_2^2}+\frac{2 \De_1^2}{n}+\frac{4 \De_3 N^2}{n N_1 N_2}.
\end{eqnarray}
By applying Chebyshev's inequality to (2.18)-(2.20), we obtain 
\begin{eqnarray}
\widetilde{\De}_0 \xrightarrow{P} \De_0,~\widetilde{\De}_1 \xrightarrow{P} \De_1.
\end{eqnarray}
Finally, from (2.16) and (2.21), we see that consistency of $\widetilde{\De}_0$ and $\widetilde{\De}_1$ imply consistency of $\Deh_0$ and $\Deh_1$. 
$\hfill\square$
\vspace{8pt}

Now by substituting the estimators of $\De_0, \De_1, a_2$ into the limiting term in Lemma 2.1. 
the consistent estimator of $e(2|1)$ is given by $\Phi((\Uh_0+c)\Vh_0^{-1/2})$, where
$\Uh_0 = - \Deh_0/2$ and $\Vh_0 = \Deh_1 + Np\hat{a}_2/(N_1N_2)$.

The following theorem is provided by the consistency of estimators $\Deh_0, \Deh_1$ and $\ah_2$,
continuous mapping theorem and dominated convergence theorem. 
\begin{thm}
Under assumptions {\rm(A1) - (A3)}
\begin{eqnarray*}
\Phi\left((\Uh_0+c)\Vh_0^{-1/2}\right)\xrightarrow{P} e(2|1)\hspace{5mm} and \hspace{5mm}
\E\left[\Phi\left((\Uh_0+c)\Vh_0^{-1/2}\right)\right]\to e(2|1).
\end{eqnarray*}
\end{thm}
\vspace{8pt}

By the results of Theorem 2.1 and Lemma 2.1, the {\bf M1}-based cut-off point for EDDR using $T(\x)$ is
 provided by 
\begin{eqnarray*}
&\ch_1 = \Vh^{1/2}_0 z_\alpha - \Uh_0,
\end{eqnarray*}
where $z_\alpha$ is the $\alpha$-percentile of $\mathcal{N}(0,1)$ and $\alpha\in (0,1)$.
\section{Asymptotic distribution of the conditional misclassification error}  
Our objective in this section is to establish the 
asymptotic distribution of $ce(2|1)$, for which we need some auxiliary notations and assumptions. 
We begin by modifying the high-dimensional asymptotic framework from Section 2 by replacing the Assumption (A3) with (B3) as follows:     
\begin{eqnarray*}
\displaystyle ({\rm B3}):\lim_{(n,p) \to \infty}\frac{\Delta_i}{n} \to 0, i=2, \cdots, 5,
~~\lim_{(n,p) \to \infty}\frac{a_i}{n} \to 0, i=3, \cdots, 6.
\end{eqnarray*}
As $ce(2|1)$ is a function of the variable set of $(U,V)$, we first obtain the joint asymptotic distribution of $(U,V)$. 
\begin{lem}
Let $\widetilde{U}=U+(N_1-N_2)p\ah_1/(N_1N_2).$ Then under assumptions $\mathrm{(A1)}$, $\mathrm{(A2)}$ and $\mathrm{(B3)}$ the following holds
\begin{eqnarray*}
\sqrt{n}\left\{
\left(
\begin{array}{c}
\widetilde{U}\\
V
\end{array}
\right)
-
\left(
\begin{array}{c}
U_0\\
V_0
\end{array}
\right)\right\}
\xrightarrow {\mathcal{D}} \mathcal{N}_2
(\boldsymbol{0},\Theta),
\end{eqnarray*}
where
\begin{eqnarray*}
\Theta=n
\left(
\begin{array}{cc}
H_U(\Delta_1,a_2)
&H_{UV}(\Delta_2,a_3)\\
H_{UV}(\Delta_2,a_3)
&H_V(\Delta_3,a_4)
\end{array}
\right),
\end{eqnarray*}
\begin{align*}
H_{UV}(\De_2,a_3)=&-\frac{2}{N_2} \Delta_2-\frac{N(N_1-N_2)pa_3}{(N_1 N_2)^2},
\end{align*}
and $\xrightarrow {\mathcal{D}}$ denotes convergence in distribution.
\end{lem}

\noindent
{\bf (Proof)}

Let $d_1$ and $d_2$ denote two non-random values which satisfy $\displaystyle  0 < \lim_{(n, p) \to \infty} |d_1| < \infty$ and $\displaystyle  
0 < \lim_{(n, p) \to \infty} |d_2| < \infty$, and introduce the statistic $Q$ which is defined as the linear combination of $\widetilde{U}$ 
and $V$. 
\begin{eqnarray*}
Q=\sqrt{n}\left\{d_1\left(\widetilde{U}+\frac{1}{2} \Delta_0\right)
+d_2\left(V-\Delta_1-\frac{Npa_2}{N_1N_2}\right)\right\}.
\end{eqnarray*}

The asymptotic normality of $Q$ would imply that the joint distribution of $\widetilde{U}$ 
and $V$ is asymptotically normal. Thus, Lemma 3.1 will be proven if we show the normal convergence 
of $Q$ under (A1), (A2) and (B3). 
We introduced the following notations
\begin{eqnarray*}
\boldsymbol{\omega}_1&=&\frac{d_1\sqrt{n}}{\sqrt{N}}\Lambda^{\frac{1}{2}}\Gamma'\bde, \\
\boldsymbol{\omega}_2&=&\frac{2d_2\sqrt{n N}}{\sqrt{N_1 N_2}} \Lambda^{3/2}\Gamma'\bde
-\frac{d_1\sqrt{n N_1}}{\sqrt{N N_2}} \Lambda^{1/2}\Gamma'\bde, \\
\Omega_3&=&\frac{d_1\sqrt{n}}{\sqrt{N_1N_2}} \Lambda, \\
\Omega_4&=&\frac{d_2 \sqrt{n} N}{N_1N_2} \Lambda^2- \frac{d_1 \sqrt{n} (N_1-N_2)}{2N_1N_2}\Lambda. 
\end{eqnarray*}
Now, since $\ah_1 -a_1= O_p(n^{-1})$ by (2.14), the statistic $Q$ can be expressed as 
\begin{eqnarray*}
Q&=&\boldsymbol{\omega}_1^\prime\textrm{\boldmath $z_1$}+\boldsymbol{\omega}_2^\prime\textrm{\boldmath $z_2$}+\textrm{\boldmath $z_1$}^\prime\Omega_3\textrm{\boldmath $z_2$}+
\textrm{\boldmath $z_2$}^\prime\Omega_4\textrm{\boldmath $z_2$}+o_p(1).
\end{eqnarray*}
Note also that 
\begin{eqnarray*}
\boldsymbol{\omega}_1'\boldsymbol{\omega}_1
&=&\frac{d_1^2n}{N}\bde'\Sigma\bde, \\
\boldsymbol{\omega}_2'\boldsymbol{\omega}_2
&=&\frac{4d_2^2nN}{N_1N_2}\bde'\Sigma^3\bde+\frac{d_1^2nN_1}{NN_2}\bde'\Sigma\bde
-\frac{4d_1d_2n}{N_2}\bde'\Sigma^2\bde,\\
\tr\Omega_3^2&=&\frac{d_1^2n}{N_1N_2}\tr\Sigma^2, \\
\tr\Omega_4^2&=&\frac{d_2^2nN^2}{N_1^2N_2^2}\tr\Sigma^4
+\frac{d_1^2n(N_1-N_2)^2}{4N_1^2N_2^2}\tr\Sigma^2
-\frac{d_1d_2n(N_1^2-N_2^2)}{N_1^2N_2^2}\tr\Sigma^3. 
\end{eqnarray*}
By combining these terms, we now obtain the asymptotic variance of $Q$ as
\begin{eqnarray*}
\sigma_Q^2=\lim_{(n,p)\to\infty}n\{d_1^2H_U(\Delta_1,a_2)-2d_1d_2H_{UV}(\Delta_2,a_3)+d_2^2H_V(\Delta_3,a_4)\} 
\end{eqnarray*}
and observe that (A1), (A2) and (B3)
\begin{eqnarray}
0<\sigma_Q^2<\infty.
\end{eqnarray}
Furthermore, the following convergence results hold
\begin{eqnarray}
\boldsymbol{\omega}_1'\Omega_3\boldsymbol{\omega}_2\to 0,
~\boldsymbol{\omega}_2'\Omega_4\boldsymbol{\omega}_2\to 0,
~\tr\Omega_3^2\Omega_4\to 0~ {\rm and} ~
\tr\Omega_4^3\to 0.
\end{eqnarray}
Now by using (3.1) and (3.2), and by applying (A.1) from Lemma A.1 (see Appendix), we obtain
\begin{eqnarray}
\frac{\boldsymbol{\omega}_1'\Omega_3\boldsymbol{\omega}_2}{\sigma_Q^3}\to 0,
~\frac{\boldsymbol{\omega}_2'\Omega_4\boldsymbol{\omega}_2}{\sigma_Q^3}\to 0,
~\frac{\tr\Omega_3^2\Omega_4}{\sigma_Q^3}\to 0~ {\rm and} ~\frac{\tr\Omega_4^3}{\sigma_Q^3}\to 0.
\label{con}
\end{eqnarray}
(3.3) in combination with Lemma A.1 show that the asymptotic normality of $Q$ holds, which completes the proof. 
$\hfill\square$

Now we are ready to state our main results on the distribution of $ce(2|1)$. Besides the distribution of the latter 
we also find the asymptotic distribution of the logit transform of $ce(2|1)$. 
Our motivation to make this particular type of transform will be clear below. 
\medskip
\begin{thm}
Let the logit transform of $ce(2|1)$ be defined by 
\begin{eqnarray*}
\ell(2|1)= \log{\frac{ce(2|1)}{1-ce(2|1)}}
\end{eqnarray*}
and let the operator $\nabla_{(u,v)} (\cdot)$ for a function $f(u,v)$ be defined as  
\begin{eqnarray*}
\nabla_{(u,v)}f(u,v) = \left ( \frac{\partial f}{\partial u} , \frac{\partial f}{\partial v} \right )^\prime. 
\end{eqnarray*}
Then in the framework $\mathrm{(A1)}$, $\mathrm{(A2)}$ and $\mathrm{(B3)}$ $ce(2|1)$ and $\ell(2|1)$ are asymptotically normal, i.e.  
\begin{eqnarray*}
&{\rm (i)}&ce(2\mid 1)\xrightarrow{\mathcal{D}} \mathcal{N} \left (e_0 ,\tau^2 \right ),\\
&{\rm (ii)}& \ell (2\mid 1) \xrightarrow{\mathcal{D}} 
\mathcal{N}\left(\ell_{0},\tau_{\ell}^2\right)
\end{eqnarray*}
with 
\begin{eqnarray*}
e_0= \Phi\left(\frac{U_0 + c}{V_0^{1/2}}\right),~\ell_{0}= \log{\frac{e_0}{1-e_0}},    ~\tau^2=\nabla_{(U_0,V_0)}' \Theta \nabla_{(U_0,V_0)}, ~
\tau^2_{\ell}=\frac{\tau^2}{(1-e_0)e_0},
\end{eqnarray*} 
where $\nabla_{(U_0,V_0)}$ is defined as 
\begin{eqnarray*}
\nabla_{(U_0,V_0)}= \left(V_0^{-1/2}\phi\left(\frac{U_0 +c }{\sqrt{V_0}} \right),
- \frac{(U_0 + c)}{2 V_0^{3/2}} \phi \left(\frac{U_0 + c}{\sqrt{V_0}} \right)\right)^\prime. 
\end{eqnarray*}
\end{thm}
\noindent
{\bf (Proof)} 
By using asymptotic normality of $(\widetilde{U},V)$ and by applying Lemma A.4 (see Appendix)  to the function
\begin{eqnarray*}
g(\widetilde{U},V)=\Phi\left(\frac{\widetilde{U}+c}{V^{1/2}}\right)
\end{eqnarray*}
it easily follows that
\begin{eqnarray*}
\nabla_{(\tilde{u},v)}g(\tilde{u},v)=\left(\frac{\partial g}{\partial \tilde{u}}, \frac{\partial g}{\partial v}\right)^\prime
= \left(v^{-1/2}\phi\left(\frac{\tilde{u} + c}{\sqrt{v}} \right),
- \frac{(\tilde{u} + c)}{2 v^{3/2}} \phi \left(\frac{\tilde{u} + c}{\sqrt{v}} \right)\right)^\prime.
\end{eqnarray*}
Then we obtain 
\begin{eqnarray*}
ce(2|1)=g(\widetilde{U},V)\xrightarrow{\mathcal{D}}\Nc(\Phi((U_0 + c)V_0^{-1/2}),\nabla_{(U_0,V_0)}' \Theta \nabla_{(U_0,V_0)}).
\end{eqnarray*}
The statement (ii) can be proven similarly. $\hfill\square$

Now we are ready to explore the determination method ${\bf M2}$ which chooses the cut-off point $c$ to 
get the desired level of confidence $1-\beta$ of a pre-specified upper bound $eu$. By the asymptotic normality of $ce(2|1)$ and $\ell(2|1)$, 
we propose to set the cut-off points for the EDDR using $T(\x)$ as
\begin{eqnarray*}
&{\rm (i)}&c_{2,1}~~s.t.~~c_{2,1}=\frac{-U_0+V_0^{1/2}z_{\gamma}}{a_1},\\
&{\rm (ii)}&c_{2,2}~~s.t.~~c_{2,2}=\frac{-U_0+V_0^{1/2}z_{\gamma_{\ell}}}{a_1},\\
\end{eqnarray*}
where
\begin{eqnarray*}
\gamma=eu-\tau z_{1-\beta},~\gamma_{\ell}=
\frac{eu}{(1-eu)\exp(\tau_{\ell} z_{1-\beta})+eu}. 
\end{eqnarray*}
\begin{remark}
If $\gamma\notin [0,1]$ then (i) is not defined. This motivates our logit trance form of $ce(2|1)$ which yields the result (ii) 
where $\gamma_\ell\in [0,1]$ always. 
\end{remark}
For practical use, the unknown parameters $\De_0$, $\De_1$, $\De_2$, $\De_3$, $a_1$, $a_2$, $a_3$ and $a_4$ in (i)-(ii) should be replaced by their consistent 
estimators. 
To ensure consistency, the asymptotic framework (A1)-(A3) is modified by replacing (A3) with
\begin{eqnarray*}
\displaystyle ({\rm B'3}):0 < \lim_{(n,p) \to \infty}a_i < \infty, ~i=3, \dots, 8, ~~0 < \lim_{(n,p) \to \infty}\De_i < \infty, ~i=2, \dots, 7.
\end{eqnarray*}

By the consistency results of Lemma A.5 and A.6 (see Appendix), obtained under the assumptions (A1), (A2) and (B$'$3), we now propose the 
{\bf M2}-based cut-off point estimator as, 
\begin{eqnarray*}
&{\rm (i)}&\hat{c}_{2,1}~~s.t.~~\hat{c}_{2,1}=\frac{-\Uh_0+\Vh_0^{1/2}z_{\hat{\gamma}}}{\ah_1},\\
&{\rm (ii)}&\hat{c}_{2,2}~~s.t.~~\hat{c}_{2,2}=\frac{-\Uh_0+\Vh_0^{1/2}z_{\hat{\gamma}_{\ell}}}{\hat{a}_1},\\
\end{eqnarray*}
where 
\begin{eqnarray*}
\hat{\gamma}=eu-\hat{\tau} z_{1-\beta},~\hat{\gamma}_\ell=\frac{eu}{(1-eu)\exp(\hat{\tau}_{\ell} z_{1-\beta})+eu}. 
\end{eqnarray*}
\begin{remark}
The problem described in Remark 3.1 remains for $\hat{\gamma}$. Therefore for practical use we recommend to replace $\hat{\gamma}$ with $\hat{\gamma_\ell}$ 
when the observed value of $\hat{\gamma}\notin [0,1]$. 
\end{remark}
\section{Simulation study}

We now turn to numerical evaluation of the asymptotic results and the suggested cut-off points. 
The goal of the simulation experiment is threefold: to investigate the finite sample behaviour of newly derived asymptotic approximations, 
to compare the performance of our approach under independence with that for dependent data with various dependence strength, and to investigate the effect of choice of the confidence level in combination with the upper bound specification.  

The data sets for each $\Pi_{g}$, $g=1,2$ are independently generated as
\begin{align}
&\boldsymbol{x}_{11},\boldsymbol{x}_{12},\ldots,
\boldsymbol{x}_{1N_1} \stackrel{i.i.d.}\sim \mathcal{N}_p(\boldsymbol{\mu}_1,\Sigma), \hspace{-1.5cm}
&\boldsymbol{x}_{21},\boldsymbol{x}_{22},\ldots,
\boldsymbol{x}_{2N_2} \stackrel{i.i.d.}\sim \mathcal{N}_p(\boldsymbol{\mu}_2,\Sigma),
\label{datasim}
\end{align}
respectively. To assess the performance  for dependent data, $\Sigma$ will be assumed to have band correlation $\Sigma=\left( \sigma_{ij} \right)$, 
\begin{align*}
\sigma_{ij} =
\begin{cases}{}
\rho^{|i-j|}, ~~~|i-j| \leq 50,\\
0,~~~~~~~~|i-j|>50,
\end{cases}
\end{align*}
with $\rho$ ranging from $0$ to $0.5$, which is chosen to fulfill the condition $(\textrm{A}2)$. To constrain the classification complexity, we set  
\begin{align*}
\Sigma^{-1/2}\boldsymbol{\mu}_1=(p)^{-1/2}(5^{1/2},5^{1/2},\ldots,5^{1/2})'~~~ {\textrm{and}} ~~~\boldsymbol{\mu}_2=(0,0,\ldots,0)',
\end{align*}
through the whole simulation experiment. 

To evaluate the effect of high-dimensionality and sample size, we let $p=64,128,256,$ 
$512,1024$ and $N_1=N_2$, $N=64,128,256$ for each choice of $\rho$. 

First, as in the previous sections, we focus without loss of generality on evaluation of $ce(2|1)$. For each triple $(p,N,\rho)$, 
we generate data according to $(\ref{datasim})$,  apply EDDR given by $T(\x)$ in (1.3) with both $\mathbf{M1}$-based cut-offs, 
$\hat{c}_1$ established in Section 2, and repeat the whole process independently $100 ~000$ times.  
As a result, we get $100 ~000$ conditional classification errors  of $T(\x)$: 
\begin{align*}
C^{(i)}=\Phi\left(\frac{U^{(i)}+(N_2^{-1}-N_1^{-1})p\hat{a}^{(i)}_1/2+\hat{c}^{(i)}_1}{\sqrt{V^{(i)}}}\right),~~~i=1,\dots,100~000,
\end{align*}
which after averaging provides {\it attained error rate} 
\begin{align*}
ae(\hat{c}_1)=\frac{1}{100~000}\sum_{i=1}^{100~000}C^{(i)}.
\end{align*}
This result, being summarized in Table 1 through Table 9, suggest that the EDDR based on $\hat{c}_1$ 
is optimally adaptive in a sense that its performance accuracy is closely approaching the actual value of the misclassification, 
$\alpha$. Stably good result is obtained when varying the dependence strength $\rho$ and the value  of the actual error $\alpha$, 
in both large sample and high-dimensional cases. 

To evaluate the performance of the $\mathbf{M2}$-based cut-offs we use the simulation setting $(\ref{datasim})$, 
with the same variety of covariance strength, a range of $\beta$ varying between $0.01$ to $0.1$ 
representing higher respective lower confidence levels, and two values of $eu$, $0.1$ and $0.2$ 
representing the upper bound on the actual misclassification  probability. We summarize the combination of the values of $1-\beta$ and $eu$ in Table 10. 
Then for each setting, 
the classification procedure by $T(\x)$ with cut-offs $\hat{c}_{2,1}$ and $\hat{c}_{2,2}$ in section 3, 
respectively. 
Proceeding with the same simulation strategy as above for each cut-off choice, we consider the {\it attained confidence level} 
\begin{align*}
acl(\hat{c}_{2,i}) = \frac{\# \left \{\Phi\left(\{U + (N_2^{-1}-N_1^{-1})\hat{a}_1/2 + \hat{c}_{2,i}\}/\sqrt{V} \right) \leq eu \right\} }{100~000},
~ i=1,2,
\end{align*}
which is obtained by averaging the observed confidence level of $ce(2|1)$ of $T(\x)$ with $\hat{c}_{2,i}$ for each, $i$, 
over $100 ~000$ independent replicates of the data generation step, estimation of parameters and classification. This result, being summarized in Table 11 through Table 28.
In most tables, the case in using $\hat{c}_{2,2}$ is better accuracy than the case in using $\hat{c}_{2,1}$, and conservative.
\section{Conclusion}

This paper contributes to the asymptotic analyses of the EDDR performance in high-dimensional data, with particular focus on determining a cut-off point to adjust the probabilities of misclassification. Two generic cut-off determination approaches, 
$\mathbf{M1}$ based on the expected error and $\mathbf{M2}$ based on the upper bound of the actual misclassification probability, $eu$ with the specified confidence level $1-\beta$, are proposed.  

To establish the cut-off by $\mathbf{M1}$, an approximation of the expected misclassification error along with its asymptotic unbiased estimator, is derived;  our result extends the approach of Anderson (1973) by considering a more general asymptotic set-up that allows $p > N$. Subsequently, the cut-off based on the main term of the asymptotic expression is suggested.

To set up the cut-off based on $\mathbf{M2}$, the asymptotic normality of the conditional misclassification error and its logit transform are established for a given $\beta$ and $eu$ in high-dimensions. 
Based on the asymptotic results, two types of cut-offs are also established. Our newly derived results extend the asymptotic consideration by McLachlan (1977) to a high-dimensional case. 

For both $\mathbf{M1}$ and $\mathbf{M2}$ approaches, the practically workable expressions of the theoretical cut-offs are established, for which we obtain consistent and asymptotic unbiased estimators of a set of unknown parameters. The validity of the new asymptotic results in a finite sample case is numerically shown by applying the cut-offs in the suggested EDDR classifier $T(\x)$ for a range of confidence levels, various strength of correlation and a 
set of $p$ and $N$ values.  

As the both suggested cut-off determination procedures demonstrate stably good accuracy in high dimensions, they can generally be recommended for practical applications in distance-based classifiers, with EDDR as special case, when it is desired to set a cut-off point to achieve a specified misclassification error.
\\

{\bf Acknowledgments.}\ \
The authors thank Professor Makoto Aoshima and Professor Yasunori Fujikoshi for extensive discussions, references and encouragements. The research of Tatjana Pavlenko is in part supported by the grant 2013-45266 VR of Sweden. The research of Takashi Seo was supported in part by Grant-in-Aid for Scientific Research (C) (23500360). The research of Masashi Hyodo is in part supported by the Stiftelsen G.S. Magnuson travel grant (2013), the Royal Swedish Academy of Sciences.  

%
%
\renewcommand{\theequation}{A.~\arabic{equation}}
\section*{Appendix}
\begin{lema}[The central limit theorem for quadratic forms]
Let $\z_1$ and $\z_2$  be independent, $\Nc_p(\boldsymbol{0},I_p)$ distributed random variables, 
$\boldsymbol{\omega}_i~(i=1,2)$ be arbitrary non-random $p$-dimensional vectors and $\Omega_i~(i=3,4)$ 
be arbitrary non-random $p\times p$ diagonal matrices. 
Define $K=\boldsymbol{\omega}_1'\z_1+\boldsymbol{\omega}_2'\z_2+\z_1^\prime\Omega_3\z_2+(\z_2'\Omega_4\z_2-\tr\Omega_4)$ 
with $\sigma_K^2=\boldsymbol{\omega}_1'\boldsymbol{\omega}_1+\boldsymbol{\omega}_2'\boldsymbol{\omega}_2
+\tr\Omega_3^2+2\tr\Omega_4^2$. If the following limiting conditions are fulfilled
\begin{eqnarray}
\frac{\boldsymbol{\omega}_1'\Omega_3\boldsymbol{\omega}_2}{\sigma_K^3}\to 0,
~\frac{\boldsymbol{\omega}_2'\Omega_4\boldsymbol{\omega}_2}{\sigma_K^3}\to 0,~
\frac{\tr\Omega_3^2\Omega_4}{\sigma_K^3}\to 0~
and ~\frac{\tr\Omega_4^3}{\sigma_K^3}\to 0,
\label{con}
\end{eqnarray}
then $K/\sigma_K\xrightarrow{\mathcal{D}}\Nc(0,1)$ as $p\to\infty$. 
\end{lema}
{\bf (Proof)}

Let $\omega_{ij} (i=1,2)$ be the {\it j-th} element of $\boldsymbol{\omega}_i, \omega_{ij} (i=3,4)$ be the {\it j-th} diagonal element of $\Omega_j$ and $z_{ij} (i=1,2)$
be the {\it j-th} element of $\boldsymbol{z}_i$. 
$K$ can be expressed as
\begin{eqnarray*}
K&=&\boldsymbol{\omega}_1'\z_1+\boldsymbol{\omega}_2'\z_2+
\z_1^\prime\Omega_3\z_2+(\z_2^\prime\Omega_4\z_2-\tr\Omega_4)\\
&=&\sum_{i=1}^p\omega_{1i}z_{1i}+\sum_{i=1}^p\omega_{2i}z_{2i}
+\sum_{i=1}^p\omega_{3i}z_{1i}z_{2i}+
\sum_{i=1}^p(\omega_{4i}z_{2i}^2-\omega_{4i}).
\end{eqnarray*}
Consider $\varepsilon_i=\omega_{1i}z_{1i}+\omega_{2i}z_{2i}+\omega_{3i}z_{1i}z_{2i}+\omega_{4i}z_{2i}^2-\omega_{4i},~
(i=1,2,\ldots,p)$ and note that $\{ \varepsilon_i \}_{i=1}^p$ is a sequence of i.i.d. random variables such that $K=\sum_{i=1}^p\varepsilon_i$ 
and the third moment of $\varepsilon_i$ is given by
\begin{eqnarray*}
\E[\varepsilon_i^3]=2(3\omega_{1i}\omega_{2i}\omega_{3i}  + 3\omega_{2i}^2\omega_{4i} + 3\omega_{3i}^2\omega_{4i} + 4\omega_{4i}^3).
\end{eqnarray*}
Then to ensure that $K/\sigma_K\xrightarrow{\mathcal{D}}\mathcal{N}(0,1)$, 
we consider the Lyapunov-based sufficient condition for the sequences $\{ \varepsilon_i \}_{i=1}^p$ which states that there exists such $\eta>0$ that
\begin{eqnarray}
\frac{\sum_{i=1}^p\E[\varepsilon_i^{2+\eta}]}{\sigma_K^{2+\eta}}\to 0~{\rm as}~p\to\infty.
\end{eqnarray}
For now, we check (A. 2) with $\eta=1$. Based on the third moment of $\varepsilon_i$, we obtain
\begin{eqnarray}
\sum_{i=1}^p\E[\varepsilon_i^3]=2( 3\boldsymbol{\omega}_1'\Omega_3\boldsymbol{\omega}_2 + 3\boldsymbol{\omega}_2'\Omega_4\boldsymbol{\omega}_2 + 3\tr\Omega_3^2\Omega_4 + 4\tr\Omega_4^3).\label{eq2}
\end{eqnarray}
From (\ref{eq2}) and the condition (\ref{con}),
\begin{eqnarray*}
\frac{\sum_{i=1}^p\E[\varepsilon_i^3]}{\sigma_K^3} \to 0
\end{eqnarray*}
as $p\to\infty$, from which the convergence $K/\sigma_K\xrightarrow{\mathcal{D}}\Nc(0,1)$ follows. $\hfill\square$
\newpage
\begin{lema}[Higher order moments of the traces of Wishart matrices]
Let $W$ be distributed as $\mathcal{W}_p(n,\Sigma)$, where $\mathcal{W}_p$ denoted Wishart distribution with freedom parameter $n$ and scale parameter $\Sigma$.
Let $A$ and $B$ denote $p \times p$ symmetric non-random matrices. 
Then the following assertions hold:
\begin{eqnarray*}
&&\mathrm{(i)}~E[(\tr A W)(\tr BW)]=n^2 \tr\Sigma A\tr\Sigma B+2 n \tr\Sigma
   A \Sigma B,\\
&&\mathrm{(ii)}~E[\tr AWBW]=(n^2+n)\tr\Sigma A \Sigma B+n \tr\Sigma A
   \tr\Sigma B,\\
&&\mathrm{(iii)}~E[(\tr A W^3)]=
(n^3+3n^2+4n) \tr\Sigma^3 A+(n^2+n) \tr\Sigma^2\tr\Sigma A\\
&&~~~~~~~~~~~~~~~~~~~~~~~
+2 n(n+1) \tr\Sigma\tr\Sigma^2 A
+n (\tr\Sigma)^2 \tr\Sigma A
,\\
&&\mathrm{(iv)}~\E[(\tr A W^2)(\tr BW^2)]=
n(n^2+n+2)(n+1)\tr\Sigma^2 A \tr\Sigma^2 B
\\
&&~~~~~~~~~~~~~~~~~~~~~~~~~~~~~~~~~~
+(\tr\Sigma)^2 \{n^2 \tr\Sigma A\tr\Sigma B+2 n \tr\Sigma A \Sigma B\}
\\
&&~~~~~~~~~~~~~~~~~~~~~~~~~~~~~~~~~~
+\tr(\Sigma) 
\{n (n^2+n+2) \tr\Sigma^2 A \tr\Sigma B
\\
&&~~~~~~~~~~~~~~~~~~~~~~~~~~~~~~~~~~
+n (n^2+n+2) \tr\Sigma A\tr\Sigma^2 B
+8 n (n+1) \tr\Sigma^2A \Sigma B\}
\\
&&~~~~~~~~~~~~~~~~~~~~~~~~~~~~~~~~~~
+\tr\Sigma^2 (2 n \tr\Sigma A \tr\Sigma B
+2 n (n+1) \tr\Sigma A \Sigma B)
\\
&&~~~~~~~~~~~~~~~~~~~~~~~~~~~~~~~~~~
+4 n (n+1)^2 \tr\Sigma^2 A \Sigma^2 B+4 n(n^2+3n+4)\tr\Sigma^3A\Sigma B
\\
&&~~~~~~~~~~~~~~~~~~~~~~~~~~~~~~~~~~
+4 n (n+1)(\tr\Sigma^3A \tr\Sigma B+\tr\Sigma A\tr\Sigma^3 B)
,\\
&&\mathrm{(v)}~\E[\tr AW^2BW^2]=
2n(n+1)^2\tr\Sigma^2 A\tr\Sigma^2 B
+n(n^2+3n+4)(n+1)\tr\Sigma^2 A \Sigma^2 B
\\
&&~~~~~~~~~~~~~~~~~~~~~~~~~~~
+(\tr\Sigma)^2(n\tr\Sigma A\tr\Sigma B+n(n+1)\tr\Sigma A \Sigma B)
\\
&&~~~~~~~~~~~~~~~~~~~~~~~~~~~
+\tr(\Sigma) \{2 n (n+1)(\tr\Sigma^2 A\tr\Sigma B+\tr\Sigma A\tr\Sigma^2 B)
\\
&&~~~~~~~~~~~~~~~~~~~~~~~~~~~
+2n(n^2+3n+4)\tr\Sigma^2 A \Sigma B\}
\\
&&~~~~~~~~~~~~~~~~~~~~~~~~~~~
+\tr(\Sigma^2)\{n(n+1) \tr\Sigma A\tr\Sigma B+n (n+3) \tr\Sigma A \Sigma B\}
\\
&&~~~~~~~~~~~~~~~~~~~~~~~~~~~
+n (n^2+3n+4)(\tr\Sigma^3 A\tr\Sigma B+\tr\Sigma A\tr\Sigma^3 B)
\\
&&~~~~~~~~~~~~~~~~~~~~~~~~~~~
+2 n (n^2+7n+8) \tr\Sigma^3 A \Sigma B, \\
&&\mathrm{(vi)}~E[(\tr A W^3)(\tr BW^3)]= (n^2\tr\Sigma  A \tr\Sigma  B+2n\tr\Sigma  A  \Sigma  B ) (\tr\Sigma)^4\\
&&~~~~~~~~~~~~~~~~~~~~~~~~~~~~~~~~~~
+ \{2 n (n^2+n+2 )(\tr\Sigma  B \tr\Sigma^2  A+\tr\Sigma  A\tr\Sigma^2 B)\\
&&~~~~~~~~~~~~~~~~~~~~~~~~~~~~~~~~~~
+16 n (n+1) \tr\Sigma^2 A  \Sigma  B \}(\tr\Sigma)^3
\\
&&~~~~~~~~~~~~~~~~~~~~~~~~~~~~~~~~~~
   + \{2 n  (n^2+n+4 ) \tr\Sigma  A
   \tr\Sigma  B\tr\Sigma^2
\\
&&~~~~~~~~~~~~~~~~~~~~~~~~~~~~~~~~~~
   +12 n (n+1) \tr\Sigma  A  \Sigma  B \tr\Sigma^2
\\
&&~~~~~~~~~~~~~~~~~~~~~~~~~~~~~~~~~~
+4 n (n+1)  (n^2+n+4 )
   \tr\Sigma^2 A \tr\Sigma^2 B
\\
&&~~~~~~~~~~~~~~~~~~~~~~~~~~~~~~~~~~
+n (n^3+3n+24n+20) (\tr\Sigma  B \tr\Sigma^3 A+\tr\Sigma  A \tr\Sigma^3  B)
\\
&&~~~~~~~~~~~~~~~~~~~~~~~~~~~~~~~~~~
+4 n (5 n^2+11n+8) \tr\Sigma^2  A  \Sigma^2 B
\\
&&~~~~~~~~~~~~~~~~~~~~~~~~~~~~~~~~~~
+24 n (n^2+3n+4) \tr\Sigma^3  A  \Sigma  B\}(\tr\Sigma)^2
\\
&&~~~~~~~~~~~~~~~~~~~~~~~~~~~~~~~~~~
+ \{2 n (n+1)  (n^2+n+10 )\\
&&~~~~~~~~~~~~~~~~~~~~~~~~~~~~~~~~~~
\times(
   \tr\Sigma  B \tr\Sigma^2 \tr\Sigma^2 A
   +\tr\Sigma  A \tr\Sigma^2 \tr\Sigma^2  B)
\\
&&~~~~~~~~~~~~~~~~~~~~~~~~~~~~~~~~~~
+2 n (n^4+4n^3+21n^2+38n+32) \tr\Sigma^3 B\tr\Sigma^2  A\\
&&~~~~~~~~~~~~~~~~~~~~~~~~~~~~~~~~~~
+16 n (n+1)
   \tr\Sigma  A\tr\Sigma  B \tr\Sigma^3
\\
&&~~~~~~~~~~~~~~~~~~~~~~~~~~~~~~~~~~
+8 n (n^2+3n+4) \tr\Sigma^3 \tr\Sigma  A \Sigma  B
\\
&&~~~~~~~~~~~~~~~~~~~~~~~~~~~~~~~~~~
+2 n (n^4+4n^3+21n^2+38n+32) \tr\Sigma^2  B \tr\Sigma^3  A
\\
&&~~~~~~~~~~~~~~~~~~~~~~~~~~~~~~~~~~
+16 n (2 n^2+5n+5) \tr\Sigma^2 \tr\Sigma^2  A  \Sigma  B
\end{eqnarray*}
\begin{eqnarray*}
&&~~~~~~~~~~~~~~~~~~~~~~~~~~~~~~~~~~
+4 n (7 n^2+19n+22) (\tr\Sigma  B \tr\Sigma^4  A+\tr\Sigma  A \tr\Sigma^4  B)
\\
&&~~~~~~~~~~~~~~~~~~~~~~~~~~~~~~~~~~
+16 n (2 n^3+9n^2+21n+16) \tr(\Sigma^3  A \Sigma^2  B)
\\
&&~~~~~~~~~~~~~~~~~~~~~~~~~~~~~~~~~~
+16 n (n^3+6n^2+21n+20) \tr\Sigma^4  A  \Sigma  B\} \tr\Sigma
\\
&&~~~~~~~~~~~~~~~~~~~~~~~~~~~~~~~~~~
+n (n+1)(n^2+n+4) \tr\Sigma  A \tr\Sigma  B (\tr\Sigma^2)^2
\\
&&~~~~~~~~~~~~~~~~~~~~~~~~~~~~~~~~~~
+4 n (5 n^2+11n+8) \tr\Sigma^2 \tr\Sigma^2 A \tr\Sigma^2 B
\\
&&~~~~~~~~~~~~~~~~~~~~~~~~~~~~~~~~~~
+4 n (3 n^2+7n+6) (\tr\Sigma  B \tr\Sigma^2 A +\tr\Sigma  A \tr\Sigma^2 B) \tr\Sigma^3
\\
&&~~~~~~~~~~~~~~~~~~~~~~~~~~~~~~~~~~
+2 n (2 n^2+5n+5) (\tr\Sigma^2)^2\tr(\Sigma A  \Sigma  B)
\\
&&~~~~~~~~~~~~~~~~~~~~~~~~~~~~~~~~~~
+n (n^4+4n^3+19n^2+36n+36)
\\
&&~~~~~~~~~~~~~~~~~~~~~~~~~~~~~~~~~~
\times(\tr\Sigma  B\tr\Sigma^3  A +\tr\Sigma  A \tr\Sigma^3  B)\tr\Sigma^2
\\
&&~~~~~~~~~~~~~~~~~~~~~~~~~~~~~~~~~~
+n (n^5+6n^4+27n^3+74n^2+156n+120) \tr\Sigma^3  A \tr\Sigma^3  B
\\
&&~~~~~~~~~~~~~~~~~~~~~~~~~~~~~~~~~~
+4 n (2 n^2+5n+5) \tr\Sigma  A\tr\Sigma  B \tr\Sigma^4
\\
&&~~~~~~~~~~~~~~~~~~~~~~~~~~~~~~~~~~
+2 n (n^3+6n^2+21n+20) \tr\Sigma  A  \Sigma  B \tr\Sigma^4
\\
&&~~~~~~~~~~~~~~~~~~~~~~~~~~~~~~~~~~
+8 n (n^3+5n^2+14n+12)\tr\Sigma^3 \tr\Sigma^2 A  \Sigma  B
\\
&&~~~~~~~~~~~~~~~~~~~~~~~~~~~~~~~~~~
+4 n (2 n^3+9n^2+21n+16)
\\
&&~~~~~~~~~~~~~~~~~~~~~~~~~~~~~~~~~~
\times(2\tr\Sigma^2 B \tr\Sigma^4 A+2\tr\Sigma^2  A \tr\Sigma^4  B+\tr\Sigma^2 \tr\Sigma^2  A  \Sigma^2  B)
\\
&&~~~~~~~~~~~~~~~~~~~~~~~~~~~~~~~~~~
+12 n (n^3+5n^2+14n+12)\\
&&~~~~~~~~~~~~~~~~~~~~~~~~~~~~~~~~~~
\times(\tr\Sigma^2\tr\Sigma^3  A  \Sigma  B+\tr\Sigma  B \tr\Sigma^5  A+\tr\Sigma  A\tr\Sigma^5  B)
\\
&&~~~~~~~~~~~~~~~~~~~~~~~~~~~~~~~~~~
+2 n (3 n^4+20n^3+77n^2+152n+132) \tr\Sigma^3 A  \Sigma^3  B
\\
&&~~~~~~~~~~~~~~~~~~~~~~~~~~~~~~~~~~
+8 n (n^4+8n^3+39n^2+80n+64) \tr\Sigma^4  A  \Sigma^2B
\\
&&~~~~~~~~~~~~~~~~~~~~~~~~~~~~~~~~~~
+4 n (n^4+10n^3+65n^2+160n+148) \tr\Sigma^5 A  \Sigma  B,\\
&&\mathrm{(vii)}~E[\tr AW^3BW^3]=\{
(n^2+n)\tr\Sigma A \Sigma B
+n\tr\Sigma A \Sigma B\}(\tr\Sigma)^4
\\
&&~~~~~~~~~~~~~~~~~~~~~~~~~~~~~
+\{4n(n+1)(\tr\Sigma B\tr\Sigma^2 A
+\tr\Sigma A\tr\Sigma^2 B)
\\
&&~~~~~~~~~~~~~~~~~~~~~~~~~~~~~
+4n(n^2+3n+4)\tr\Sigma^2 A \Sigma B\} (\tr\Sigma)^3
\\
&&~~~~~~~~~~~~~~~~~~~~~~~~~~~~~
+\{6n(n+1)\tr\Sigma A\tr\Sigma B \tr\Sigma^2
\\
&&~~~~~~~~~~~~~~~~~~~~~~~~~~~~~
+2n(n^2+4n+7)\tr\Sigma A \Sigma B\tr\Sigma^2
\\
&&~~~~~~~~~~~~~~~~~~~~~~~~~~~~~
+2n(5n^2+11n+8)\tr\Sigma^2 A \tr\Sigma^2 B
\\
&&~~~~~~~~~~~~~~~~~~~~~~~~~~~~~
+6n(n^2+3n+4)(\tr\Sigma B\tr\Sigma^3 A+\tr\Sigma A\tr\Sigma^3 B)
\\
&&~~~~~~~~~~~~~~~~~~~~~~~~~~~~~
+2n(2 n^3+9n^2+21n+16) \tr\Sigma^2 A \Sigma^2 B
\\
&&~~~~~~~~~~~~~~~~~~~~~~~~~~~~~
+2 n(n^3+9n^2+42n+44) \tr\Sigma^3 A \Sigma B\} (\tr\Sigma)^2
\\
&&~~~~~~~~~~~~~~~~~~~~~~~~~~~~~
+\{4 n (2 n^2+5n+5)(\tr\Sigma B\tr\Sigma^2 A+\tr\Sigma A \tr\Sigma^2 B)\tr\Sigma^2 
\\
&&~~~~~~~~~~~~~~~~~~~~~~~~~~~~~
+4n((2 n^3+9n^2+21n+16)(\tr\Sigma^3 B\tr\Sigma^2 A+\tr\Sigma^2 B\tr\Sigma^3 A)
\\
&&~~~~~~~~~~~~~~~~~~~~~~~~~~~~~
+4 n (n^2+3n+4) \tr\Sigma A\tr\Sigma B\tr\Sigma^3
\\
&&~~~~~~~~~~~~~~~~~~~~~~~~~~~~~
+4 n (n^2+7n+8) \tr\Sigma^3\tr\Sigma A \Sigma B
\\
&&~~~~~~~~~~~~~~~~~~~~~~~~~~~~~
+4 n (n^3+6n^2+21n+20)\\
&&~~~~~~~~~~~~~~~~~~~~~~~~~~~~~
\times(\tr\Sigma^2\tr\Sigma^2 A\Sigma B+\tr\Sigma B\tr\Sigma^4 A+\tr\Sigma A\tr\Sigma^4 B)
\\
&&~~~~~~~~~~~~~~~~~~~~~~~~~~~~~
+4 n (n^4+8n^3+39n^2+80n+64) \tr\Sigma^3 A \Sigma^2B
\\
&&~~~~~~~~~~~~~~~~~~~~~~~~~~~~~
+8 n (n^3+13n^2+40n+42) \tr\Sigma^4 A \Sigma B\}\tr\Sigma
\\
&&~~~~~~~~~~~~~~~~~~~~~~~~~~~~~
+n (2 n^2+5n+5) \tr\Sigma A\tr\Sigma B (\tr\Sigma^2)^2
\end{eqnarray*}
\begin{eqnarray*}
&&~~~~~~~~~~~~~~~~~~~~~~~~~~~~~
+2 n (2 n^3+9n^2+21n+16) \tr\Sigma^2\tr\Sigma^2 A\tr\Sigma^2B
\\
&&~~~~~~~~~~~~~~~~~~~~~~~~~~~~~
+n(n^3+5n^2+14n+12)\\
&&~~~~~~~~~~~~~~~~~~~~~~~~~~~~~
\times\{2(\tr\Sigma B\tr\Sigma^2 A+\tr\Sigma A\tr\Sigma^2 B)\tr\Sigma^3\\
&&~~~~~~~~~~~~~~~~~~~~~~~~~~~~~
+3(\tr\Sigma B\tr\Sigma^2\tr\Sigma^3 A+\tr\Sigma A\tr\Sigma^2\tr\Sigma^3 B)\}
\\
&&~~~~~~~~~~~~~~~~~~~~~~~~~~~~~
+n (n^3+4n^2+10n+9) (\tr\Sigma^2)^2 \tr\Sigma A \Sigma B
\\
&&~~~~~~~~~~~~~~~~~~~~~~~~~~~~~
+n (3 n^4+20n^3+77n^2+152n+132) \tr\Sigma^3 A \tr\Sigma^3 B
\\
&&~~~~~~~~~~~~~~~~~~~~~~~~~~~~~
+n (n^3+6n^2+21n+20) \tr\Sigma A \tr\Sigma B \tr\Sigma^4
\\
&&~~~~~~~~~~~~~~~~~~~~~~~~~~~~~
+n (n^3+14n^2+41n+40) \tr\Sigma A \Sigma B\tr\Sigma^4
\\
&&~~~~~~~~~~~~~~~~~~~~~~~~~~~~~
+4 n (n^3+11n^2+28n+24) \tr\Sigma^3\tr\Sigma^2 A \Sigma B
\\
&&~~~~~~~~~~~~~~~~~~~~~~~~~~~~~
+2 n (n^4+8n^3+39n^2+80n+64)\\
&&~~~~~~~~~~~~~~~~~~~~~~~~~~~~~
\times(\tr\Sigma^2 B\tr\Sigma^4 A+\tr\Sigma^2 A\tr\Sigma^4 B)
\\
&&~~~~~~~~~~~~~~~~~~~~~~~~~~~~~
+2 n (2n^3+19n^2+43n+32) \tr\Sigma^2\tr\Sigma^2 A \Sigma^2 B
\\
&&~~~~~~~~~~~~~~~~~~~~~~~~~~~~~
+2 n (n^5+7n^4+34n^3+78n^2+72) \tr\Sigma^2\tr\Sigma^3 A \Sigma B
\\
&&~~~~~~~~~~~~~~~~~~~~~~~~~~~~~
+n (n^4+10n^3+65n^2+160n+148)\\
&&~~~~~~~~~~~~~~~~~~~~~~~~~~~~~
\times (\tr\Sigma B\tr\Sigma^5 A+\tr\Sigma A\tr\Sigma^5 B)
\\
&&~~~~~~~~~~~~~~~~~~~~~~~~~~~~~
+n (n^5+9n^4+47n^3+151n^2+308n+252) \tr\Sigma^3 A \Sigma^3 B
\\
&&~~~~~~~~~~~~~~~~~~~~~~~~~~~~~
+4 n (n^4+16n^3+75n^2+164n+128) \tr\Sigma^4 A \Sigma^2 B
\\
&&~~~~~~~~~~~~~~~~~~~~~~~~~~~~~
+2 n (n^4+22n^3+125n^2+328n+292) \tr\Sigma^5 A \Sigma B.\\
\end{eqnarray*}
\end{lema}
\vspace{-6mm}
{\bf (Proof)} 
The proof of assertions (i)-(vii) follows directly by applying 
the technique derived in Lemma A.2, in Hyodo et al. (2012). $\hfill\square$

\begin{lema}[Moments of quadratic form]
Let $\boldsymbol{x}$ be distributed $\mathcal{N}_p(\boldsymbol{0},I_p)$. 
Then the following assertions hold:
\begin{eqnarray*}
&{\rm (i)}&~\E[\x'A\x]=\tr A,\\
&{\rm (ii)}&~\E[\x' A\x\x' B\x]=
2\tr AB+\tr A\tr B,
\end{eqnarray*}
where $A$ and $B$ are $p\times p$ non-random symmetric matrices. 
\end{lema}
{\bf (Proof)} See, Gupta and Nagar (1999). 
\begin{lema}[Multivariate Delta Method]
Suppose that $\y_n=(y_{n1},\ldots,y_{nk})'$ is a 
sequence of the random vectors such that
\begin{eqnarray*}
\sqrt{n}(\y_n-\bmu)\xrightarrow{\mathcal{D}}\Nc_k(\boldsymbol{0},\Theta)~{as}~n\to\infty,
\end{eqnarray*}
where $\bmu=(\mu_1,\cdots,\mu_k)'$ is the asymptotic mean vector and $\Theta$ is the asymptotic covariance matrix 
which is assumed to be positive definite. 
Let $g:\mathbb{R}^k\to \mathbb{R}$ and g is continuously differentiable. Let
\begin{eqnarray*}
\nabla g(\y)=\left(\frac{\partial g}{\partial y_1},\cdots,\frac{\partial g}{\partial y_k}\right)'.
\end{eqnarray*}
Let $\nabla_\mu$ denote $\nabla g(\y)$ evaluated at $\y=\bmu$ and assume that the elements of 
$\nabla_\mu$ are nonzero. Then it holds that
\begin{eqnarray*}
\sqrt{n}(g(\y_n)-g(\bmu))\to\Nc(0,\nabla_\mu'\Theta \nabla_\mu)~{as}~n\to\infty.
\end{eqnarray*}
\end{lema}
{\bf (Proof)} See, Rao (1973). 
\begin{lema} [The consistent estimators of $a_3$ and $a_4$]
The consistent estimators of $a_3$ and $a_4$ are
\begin{eqnarray*}
\hat{a}_3&=&\frac{n^2}{(n+4)(n+2)(n-1)(n-2)p} \{ n^2\tr{S^3} - 3n\tr{S^2}\tr{S}+ 2(\tr{S})^3 \},\\
\hat{a}_4&=&\frac{1}{p} \{ b_1 \tr{S^4}+ b_2 \tr{S^3}\tr{S}+ b_3 (\tr{S^2})^2+ b_4 (\tr{S})^2\tr{S^2}+b_5 (\tr{S})^4 \},
\end{eqnarray*}
where
\begin{eqnarray*}
b_1&=&\frac{n^5(n^2 +n +2)}{(n+6)(n+4)(n+2)(n+1)(n-1)(n-2)(n-3)},\\
b_2&=&-\frac{4n^4(n^2+n+2)}{(n+6)(n+4)(n+2)(n+1)(n-1)(n-2)(n-3)},\\
b_3&=&-\frac{n^4(2n^2+3n-6)}{(n+6)(n+4)(n+2)(n+1)(n-1)(n-2)(n-3)},\\
b_4&=&\frac{2n^4(5n+6)}{(n+6)(n+4)(n+2)(n+1)(n-1)(n-2)(n-3)},\\
b_5&=&-\frac{n^3(5n+6)}{(n+6)(n+4)(n+2)(n+1)(n-1)(n-2)(n-3)}.
\end{eqnarray*}
\end{lema}
{\bf (Proof)} See, Hyodo et al. (2012).
\begin{lema} [The consistent estimators of $\De_2$ and $\De_3$]
The consistent estimators of $\De_2$ and $\De_3$ are
\begin{eqnarray*}
\widehat{\Delta}_2&=&\left(1+\frac{1}{n}\right)^{-1}
\left\{
\bdeh^\prime S^2\bdeh
-\frac{p}{n}\ah_1\Deh_1
-\frac{Np}{N_1N_2}
\left(\frac{n+1}{n}\ah_3+\frac{p}{n}\ah_1\ah_2\right)\right\}
,\\
\widehat{\Delta}_3&=&
\left(
\frac{n(n+3)+4}{n^2}
\right)^{-1}
\left\{\bdeh^\prime S^3\bdeh
-\frac{(n+1)p}{n^2}\ah_2\Deh_1^2
-\frac{2(n+1)p}{n^2}\ah_1\Deh_2^2
-\frac{p^2}{n^2}\ah_1^2\Deh_1^2
\right.\\
& &
-\left.\frac{Np}{N_1N_2}
\left(
\frac{n(n+3)+4}{n^2}\ah_4
+\frac{(n+1)p}{n^2}\ah_2^2
+\frac{2(n+1)p}{n^2}\ah_1\ah_3
+\frac{p^2}{n^2}\ah_1^2\ah_2
\right)\right\}.
\end{eqnarray*}
\end{lema}
{\bf (Proof)}We consider following random variables
\begin{eqnarray*}
\widetilde{\Delta}_2&=&\left(1+\frac{1}{n}\right)^{-1}
\left\{
\bdeh^\prime S^2\bdeh
-\frac{p}{n}a_1\De_1
-\frac{Np}{N_1N_2}
\left(\frac{n+1}{n}a_3+\frac{p}{n}a_1a_2\right)\right\}
,\\
\widetilde{\Delta}_3&=&
\left(
\frac{n(n+3)+4}{n^2}
\right)^{-1}
\left\{\bdeh^\prime S^3\bdeh
-\frac{(n+1)p}{n^2}a_2\De_1
-\frac{2(n+1)p}{n^2}a_1\De_2
-\frac{p^2}{n^2}a_1^2\De_1
\right.\\
& &
-\left.\frac{Np}{N_1N_2}
\left(
\frac{n(n+3)+4}{n^2}a_4
+\frac{(n+1)p}{n^2}a_2^2
+\frac{2(n+1)p}{n^2}a_1a_3
+\frac{p^2}{n^2}a_1^2a_2
\right)\right\}.
\end{eqnarray*}
Then the conditional moments of $\bdeh^\prime S^i\bdeh$ are given by
\begin{eqnarray}
E[\bdeh^\prime S^i\bdeh|S]&=&
\bde'S^i\bde
+\frac{N}{N_1N_2}\tr\Sigma S^i,\\
E[(\bdeh^\prime S^i\bdeh)^2|S]
&=&(\bde'S^i\bde)^2+
\frac{4N}{N_1N_2}\bde'S^i\Sigma S^i\bde
+\frac{N^2}{N_1^2N_2^2}\{2(\tr\Sigma S^i\Sigma S^i)+(\tr\Sigma S^i)^2\}\nonumber\\
& &+\frac{2N}{N_1N_2}\tr\Sigma S^i\bde'S^i\bde.
\end{eqnarray}
By using (A.4), (A.5) and Lemma A.2, we can calculate the expectations 
of $\widetilde{\Delta}_2,~\widetilde{\Delta}_3$ and these variances
\begin{eqnarray*}
E[\widetilde{\Delta}_2]&=&\bde^\prime \Sigma^2 \bde,~
E[\widetilde{\Delta}_3]=\bde^\prime \Sigma^3 \bde,\\
{\rm Var}[\widetilde{\Delta}_2]&=&\frac{n^2}{(n+1)^2} \left \{ \frac{4N}{N_1 N_2} \De_5 + \frac{8 N p}{n N_1 N_2} \De_4 a_1 + \frac{4}{n} \De_1 \De_3 + \frac{4 N p}{n N_1 N_2} \De_3 a_2  \right. \\
                                  && + \frac{4 N p^2}{n^2 N_1 N_2} \De_3 a_1^2 + \frac{4}{n} \De_2^2 + \frac{8 p}{n^2} \De_1 \De_2 a_1 + \frac{8 N p}{n N_1 N_2} \De_2 a_3 + \frac{8 N p^2}{n^2 N_1 N_2} \De_2 a_1 a_2  \\
                                  && + \frac{2 p}{n^2} \De_1^2 a_2 + \frac{2 p^2}{n^3} \De_1^2 a_1^2 + \frac{4 N p}{n N_1 N_2} \De_1 a_4 + \frac{8 N p^2}{n^2 N_1 N_2} \De_1 a_1 a_3 + \frac{4 N p^2}{n^2 N_1 N_2} \De_1 a_2^2 \\
                                  && + \frac{4 N p^3}{n^3 N_1 N_2} \De_1 a_1^2 a_2 + \frac{2 N^2 p}{N_1^2 N_2^2} a_6 + \frac{4 N^2 p^2}{n N_1^2 N_2^2} a_1 a_5 + \frac{4 N^2 p^2}{n N_1^2 N_2^2}  a_2 a_4 \\
                                  && + \frac{2 N^2 p^3}{n^2 N_1^2 N_2^2} a_1^2 a_4 + \frac{4 N^2 p^2}{n N_1^2 N_2^2} a_3^2 + \frac{8 N^2 p^3}{n^2 N_1^2 N_2^2} a_1 a_2 a_3 + \frac{2 N^2 p^3}{n^2 N_1^2 N_2^2} a_2^3 \\
                                  && \left. + \frac{2 N^2 p^4}{n^3 N_1^2 N_2^2} a_1^2 a_2^2  \right\}+ o(n^{-1}), \\
{\rm Var}[\widetilde{\Delta}_3]&=&\frac{n^4}{(n^2 + 3 n + 4)^2} \left \{ \frac{4 N}{N_1 N_2} \De_7 + \frac{16 N p}{n N_1 N_2} \De_6 a_1 + \frac{4}{n} \De_1 \De_5 + \frac{12 N p}{n N_1 N_2} \De_5 a_2 \right.\\
                                  && + \frac{24 N p^2}{n^2 N_1 N_2} \De_5 a_1^2 + \frac{8}{n} \De_2 \De_4 + \frac{16 p}{n^2} \De_1 \De_4 a_1 + \frac{8 N p}{n N_1 N_2} \De_4 a_3 \\
                                  && + \frac{32 N p^2}{n^2 N_1 N_2} \De_4 a_1 a_2 + \frac{16 N p^3}{n^3 N_1 N_2} \De_4 a_1^3 + \frac{6}{n} \De_3^2 + \frac{32 p}{n^2} \De_2 \De_3 a_1 \\
                                  && + \frac{12 p}{n^2} \De_1 \De_3 a_2 + \frac{24 p^2}{n^3} \De_1 \De_3 a_1^2 + \frac{12 N p}{n N_1 N_2} \De_3 a_4 + \frac{32 N p^2}{n^2 N_1 N_2} \De_3 a_1 a_3 \\
                                  && + \frac{16 N p^2}{n^2 N_1 N_2} \De_3 a_2^2 + \frac{32 N p^3}{n^3 N_1 N_2} \De_3 a_1^2 a_2 + \frac{4 N p^4}{n^4 N_1 N_2} \De_3 a_1^4 + \frac{8 p}{n^2} \De_2^2 a_2 \\
                                  && + \frac{20 p^2}{n^3} \De_2^2 a_1^2+ \frac{8 p}{n^2} \De_1 \De_2 a_3 + \frac{32 p^2}{n^3} \De_1 \De_2 a_1 a_2 + \frac{16 p^3}{n^4} \De_1 \De_2 a_1^3 \\
                                  && + \frac{8 N p}{n N_1 N_2} \De_2 a_5 + \frac{32 N p^2}{n^2 N_1 N_2} \De_2 a_1 a_4 + \frac{24 N p^2}{n^2 N_1 N_2} \De_2 a_2 a_3 + \frac{40 N p^3}{n^3 N_1 N_2} \De_2 a_1^2 a_3 \\
                                  && + \frac{32 N p^3}{n^3 N_1 N_2} \De_2 a_1 a_2^2 + \frac{16 N p^4}{n^4 N_1 N_2} \De_2 a_1^3 a_2 + \frac{2 p}{n^2} \De_1^2 a_4 + \frac{8 p^2}{n^3} \De_1^2 a_1 a_3 
\end{eqnarray*}
\begin{eqnarray*}
                                  && + \frac{4 p^2}{n^3} \De_1^2 a_2^2 + \frac{12 p^3}{n^4} \De_1^2 a_1^2 a_2 + \frac{2 p^4}{n^5} \De_1^2 a_1^4 + \frac{4 N p}{n N_1 N_2} \De_1 a_6 + \frac{16 N p^2}{n^2 N_1 N_2} \De_1 a_1 a_5 \\
                                  && + \frac{16 N p^2}{n^2 N_1 N_2} \De_1 a_2 a_4 +  \frac{24 N p^3}{n^3 N_1 N_2} \De_1 a_1^2 a_4 + \frac{8 N p^2}{n^2 N_1 N_2} \De_1 a_3^2 \\
                                  && + \frac{48 N p^3}{n^3 N_1 N_2} \De_1 a_1 a_2 a_3 + \frac{16 N p^4}{n^4 N_1 N_2} \De_1 a_1^3 a_3 + \frac{8 N p^3}{n^3 N_1 N_2} \De_1 a_2^3 \\
                                  && + \frac{24 N p^4}{n^4 N_1 N_2} \De_1 a_1^2 a_2^2 + \frac{4 N p^5}{n^5 N_1 N_2} \De_1 a_1^4 a_2 + \frac{2 N^2 p}{N_1^2 N_2^2} a_8 + \frac{8 N^2 p^2}{n N_1^2 N_2^2} a_1 a_7 \\
                                  && + \frac{8 N^2 p^2}{n N_1^2 N_2^2} a_2 a_6 + \frac{12 N^2 p^3}{n^2 N_1^2 N_2^2} a_1^2 a_6 + \frac{8 N^2 p^2}{n N_1^2 N_2^2} a_3 a_5 + \frac{24 N^2 p^3}{n^2 N_1^2 N_2^2} a_1 a_2 a_5 \\
                                  && + \frac{8 N^2 p^4}{n^3 N_1^2 N_2^2} a_1^3 a_5 + \frac{6 N^2 p^2}{n N_1^2 N_2^2} a_4^2 + \frac{32 N^2 p^3}{n^2 N_1^2 N_2^2} a_1 a_3 a_4 + \frac{16 N^2 p^3}{n^2 N_1^2 N_2^2} a_2^2 a_4 \\
                                  && + \frac{28 N^2 p^4}{n^3 N_1^2 N_2^2} a_1^2 a_2 a_4 + \frac{2 N^2 p^5}{n^4 N_1^2 N_2^2} a_1^4 a_4 + \frac{16 N^2 p^3}{n^2 N_1^2 N_2^2} a_2 a_3^2 \\
                                  && + \frac{20 N^2 p^4}{n^3 N_1^2 N_2^2} a_1^2 a_3^2 + \frac{40 N^2 p^4}{n^3 N_1^2 N_2^2} a_1 a_2^2 a_3 + \frac{16 N^2 p^5}{n^4 N_1^2 N_2^2} + \frac{4 N^2 p^4}{n^3 N_1^2 N_2^2} a_2^4 \\
                                  &&\left. + \frac{12 N^2 p^5}{n^4 N_1^2 N_2^2} a_1^2 a_2^3 +  \frac{2 N^2 p^6}{n^5 N_1^2 N_2^2} a_1^4 a_2^2 \right\} + o(n^{-1}).
\end{eqnarray*}
Using the {\it Chebyshev's inequality}, We get $\widetilde{\Delta}_i\xrightarrow{P}\De_i,~i=2,3$.

Replacing the unknown values in $\widetilde{\De}_2$ with their consistent estimator of $\De_2$, we have
\begin{eqnarray*}
\widehat{\Delta}_2&=&\left(1+\frac{1}{n}\right)^{-1}
\left\{
\bdeh^\prime S^2\bdeh
-\frac{p}{n}\ah_1\Deh_1
-\frac{Np}{N_1N_2}
\left(\frac{n+1}{n}\ah_3+\frac{p}{n}\ah_1\ah_2\right)\right\}.
\end{eqnarray*}
Using consistency of $\ah_i$, i=1,2,3 and $\Deh_1$, we can prove the consistency of $\Deh_2$.
The consistency of $\Deh_3$ can be proven similarly. $\hfill\square$

\newpage

\begin{table}[bp]
\caption{$\Sigma=I_p,~~~\alpha=0.1$}
\begin{center}
\begin{tabular}{|r|r|r|r|r|r|} \hline
$N \backslash p$ & 64 & 128 & 256 & 512 & 1024 \\\hline
64 & 0.101258 & 0.101550 & 0.101822 & 0.101591 & 0.101701 \\\hline
128 & 0.100356 & 0.100635 & 0.100767 & 0.100808 & 0.100928 \\\hline
256 & 0.100135 & 0.100175 & 0.100320 & 0.100349 & 0.100449 \\\hline
\end{tabular}
\end{center}
\label{}
\end{table}

\begin{table}[bp]
\caption{$\Sigma=I_p,~~~\alpha=0.2$}
\begin{center}
\begin{tabular}{|r|r|r|r|r|r|} \hline
$N \backslash p$ & 64  & 128  & 256  & 512 & 1024 \\\hline
64 & 0.201793  & 0.202039  & 0.201824  & 0.201865  & 0.201870  \\\hline
128 & 0.200916  & 0.201113  & 0.201099  & 0.200896  & 0.200811  \\\hline
256 & 0.200412  & 0.200442  & 0.200533  & 0.200539  & 0.200460  \\\hline
\end{tabular}
\end{center}
\label{}
\end{table}

\begin{table}[bp]
\caption{$\Sigma=I_p,~~~\alpha=0.3$}
\begin{center}
\begin{tabular}{|r|r|r|r|r|r|} \hline
$N \backslash p$ & 64 & 128 & 256 & 512 & 1024 \\\hline
64 & 0.302290 & 0.302157 & 0.301731 & 0.301633 & 0.301501 \\\hline
128 & 0.301158 & 0.301198 & 0.301026 & 0.300898 & 0.300719 \\\hline
256 & 0.300731 & 0.300555 & 0.300589 & 0.300503 & 0.300396 \\\hline
\end{tabular}
\end{center}
\label{}
\end{table}

\begin{table}[bp]
\caption{$\rho=0.2,~~~\alpha=0.1$}
\begin{center}
\begin{tabular}{|r|r|r|r|r|r|} \hline
$N \backslash p$ & 64 & 128 & 256 & 512 & 1024 \\\hline
64 & 0.100874 & 0.101402 & 0.101553 & 0.101538 & 0.101828 \\\hline
128 & 0.100087 & 0.100525 & 0.100703 & 0.100816 & 0.100883 \\\hline
256 & 0.099999 & 0.100064 & 0.100243 & 0.100349 & 0.100432 \\\hline
\end{tabular}
\end{center}
\label{}
\end{table}

\begin{table}[bp]
\caption{$\rho=0.2,~~~\alpha=0.2$}
\begin{center}
\begin{tabular}{|r|r|r|r|r|r|} \hline
$N \backslash p$ & 64 & 128 & 256 & 512 & 1024 \\\hline
64 & 0.201434  & 0.201462  & 0.201832  & 0.201741  & 0.201876  \\\hline
128 & 0.200718  & 0.201274  & 0.200788  & 0.200826  & 0.200614  \\\hline
256 & 0.200245  & 0.200183  & 0.200339  & 0.200457  & 0.200174  \\\hline
\end{tabular}
\end{center}
\label{}
\end{table}

\begin{table}[bp]
\caption{$\rho=0.2,~~~\alpha=0.3$}
\begin{center}
\begin{tabular}{|r|r|r|r|r|r|} \hline
$N \backslash p$ & 64 & 128 & 256 & 512 & 1024 \\\hline
64 & 0.302047 & 0.301750 & 0.301613 & 0.301398 & 0.301309 \\\hline
128 & 0.301085 & 0.300956 & 0.300938 & 0.300820 & 0.300959 \\\hline
256 & 0.300625 & 0.300506 & 0.300451 & 0.300399 & 0.300528 \\\hline
\end{tabular}
\end{center}
\label{}
\end{table}

\begin{table}[bp]
\caption{$\rho=0.5,~~~\alpha=0.1$}
\begin{center}
\begin{tabular}{|r|r|r|r|r|r|} \hline
$N \backslash p$ & 64 & 128 & 256 & 512 & 1024 \\\hline
64 & 0.098922 & 0.100218 & 0.100929 & 0.101420 & 0.101548 \\\hline
128 & 0.098741 & 0.099486 & 0.100154 & 0.100438 & 0.100722 \\\hline
256 & 0.099124 & 0.099461 & 0.099723 & 0.099997 & 0.100195 \\\hline
\end{tabular}
\end{center}
\label{}
\end{table}

\begin{table}[bp]
\caption{$\rho=0.5,~~~\alpha=0.2$}
\begin{center}
\begin{tabular}{|r|r|r|r|r|r|} \hline
$N \backslash p$ & 64 & 128 & 256 & 512 & 1024 \\\hline
64 & 0.202088 & 0.200184 & 0.200799 & 0.201209 & 0.201662 \\\hline
128 & 0.200979 & 0.199646 & 0.200136 & 0.200289 & 0.200633 \\\hline
256 & 0.200343 & 0.199499 & 0.199728 & 0.199936 & 0.200197 \\\hline
\end{tabular}
\end{center}
\label{}
\end{table}

\begin{table}[bp]
\caption{$\rho=0.5,~~~\alpha=0.3$}
\begin{center}
\begin{tabular}{|r|r|r|r|r|r|} \hline
$N \backslash p$ & 64 & 128 & 256 & 512 & 1024 \\\hline
64 & 0.300322 & 0.300724 & 0.300709 & 0.301292 & 0.301143 \\\hline
128 & 0.300083 & 0.300311 & 0.300173 & 0.300521 & 0.300517 \\\hline
256 & 0.299821 & 0.299983 & 0.300022 & 0.300112 & 0.300352 \\\hline
\end{tabular}
\end{center}
\label{}
\end{table}

\begin{table}[bp]
\caption{The selected parameters for cases 1-6}
\begin{center}
\begin{tabular}{|r|r|r|r|r|r|r|} \hline
 & case1 & case2 & case3 & case4 & case5 & case6 \\\hline
1-$\beta$ & 0.90  & 0.95  & 0.99  & 0.90  & 0.95  & 0.99  \\\hline
$eu$ & 0.20  & 0.20  & 0.20  & 0.10  & 0.10  & 0.10  \\\hline
\end{tabular}
\end{center}
\label{}
\end{table}

\begin{table}[bp]
\caption{$\Sigma=I_p,~\beta=0.10,~eu=0.20$}
\begin{center}
\begin{tabular}{|rr|r|r|r|r|r|r|} \hline
 &  & $p$ & 64 & 128 & 256 & 512 & 1024 \\\hline
$N$ & 64 & $acl(\hat{c}_{2,1})$ & 0.884 & 0.885 & 0.882 & 0.886 & 0.883 \\
 &  & $acl(\hat{c}_{2,2})$ & 0.906 & 0.907 & 0.904 & 0.907 & 0.905 \\\hline
 & 128 & $acl(\hat{c}_{2,1})$ & 0.891 & 0.889 & 0.888 & 0.888 & 0.889 \\
 &  & $acl(\hat{c}_{2,2})$ & 0.906 & 0.904 & 0.904 & 0.904 & 0.904 \\\hline
 & 256 & $acl(\hat{c}_{2,1})$ & 0.893 & 0.895 & 0.892 & 0.893 & 0.893 \\
 &  & $acl(\hat{c}_{2,2})$ & 0.904 & 0.905 & 0.903 & 0.904 & 0.903 \\\hline
\end{tabular}
\end{center}
\label{}
\end{table}

\begin{table}[bp]
\caption{$\Sigma=I_p,~\beta=0.05,~eu=0.20$}
\begin{center}
\begin{tabular}{|rr|r|r|r|r|r|r|} \hline
 &  & $p$ & 64 & 128 & 256 & 512 & 1024 \\\hline
$N$ & 64 & $acl(\hat{c}_{2,1})$ & 0.934 & 0.931 & 0.935 & 0.934 & 0.933 \\
 &  & $acl(\hat{c}_{2,2})$ & 0.956 & 0.954 & 0.956 & 0.956 & 0.955 \\\hline
 & 128 & $acl(\hat{c}_{2,1})$ & 0.939 & 0.938 & 0.938 & 0.938 & 0.938 \\
 &  & $acl(\hat{c}_{2,2})$ & 0.954 & 0.953 & 0.954 & 0.953 & 0.953 \\\hline
 & 256 & $acl(\hat{c}_{2,1})$ & 0.944 & 0.943 & 0.942 & 0.941 & 0.942 \\
 &  & $acl(\hat{c}_{2,2})$ & 0.954 & 0.953 & 0.952 & 0.952 & 0.953 \\\hline
\end{tabular}
\end{center}
\label{}
\end{table}

\begin{table}[bp]
\caption{$\Sigma=I_p,~\beta=0.01,~eu=0.20$}
\begin{center}
\begin{tabular}{|rr|r|r|r|r|r|r|} \hline
 &  & $p$ & 64 & 128 & 256 & 512 & 1024 \\\hline
$N$ & 64 & $acl(\hat{c}_{2,1})$ & 0.980 & 0.979 & 0.979 & 0.979 & 0.979 \\
 &  & $acl(\hat{c}_{2,2})$ & 0.992 & 0.993 & 0.993 & 0.992 & 0.992 \\\hline
 & 128 & $acl(\hat{c}_{2,1})$ & 0.983 & 0.983 & 0.983 & 0.983 & 0.982 \\
 &  & $acl(\hat{c}_{2,2})$ & 0.992 & 0.991 & 0.992 & 0.992 & 0.992 \\\hline
 & 256 & $acl(\hat{c}_{2,1})$ & 0.986 & 0.985 & 0.985 & 0.986 & 0.985 \\
 &  & $acl(\hat{c}_{2,2})$ & 0.991 & 0.991 & 0.991 & 0.991 & 0.991 \\\hline
\end{tabular}
\end{center}
\label{}
\end{table}

\begin{table}[bp]
\caption{$\Sigma=I_p,~\beta=0.10,~eu=0.10$}
\begin{center}
\begin{tabular}{|rr|r|r|r|r|r|r|} \hline
 &  & $p$ & 64 & 128 & 256 & 512 & 1024 \\\hline
$N$ & 64 & $acl(\hat{c}_{2,1})$ & 0.882 & 0.881 & 0.879 & 0.877 & 0.876 \\
 &  & $acl(\hat{c}_{2,2})$ & 0.913 & 0.911 & 0.909 & 0.907 & 0.906 \\\hline
 & 128 & $acl(\hat{c}_{2,1})$ & 0.889 & 0.887 & 0.884 & 0.884 & 0.885 \\
 &  & $acl(\hat{c}_{2,2})$ & 0.910 & 0.909 & 0.907 & 0.905 & 0.907 \\\hline
 & 256 & $acl(\hat{c}_{2,1})$ & 0.894 & 0.893 & 0.890 & 0.890 & 0.888 \\
 &  & $acl(\hat{c}_{2,2})$ & 0.909 & 0.907 & 0.905 & 0.904 & 0.903 \\\hline
\end{tabular}
\end{center}
\label{}
\end{table}

\begin{table}[bp]
\caption{$\Sigma=I_p,~\beta=0.05,~eu=0.10$}
\begin{center}
\begin{tabular}{|rr|r|r|r|r|r|r|} \hline
 &  & $p$ & 64 & 128 & 256 & 512 & 1024 \\\hline
$N$ & 64 & $acl(\hat{c}_{2,1})$ & 0.930 & 0.927 & 0.927 & 0.927 & 0.924 \\
 &  & $acl(\hat{c}_{2,2})$ & 0.960 & 0.958 & 0.958 & 0.958 & 0.956 \\\hline
 & 128 & $acl(\hat{c}_{2,1})$ & 0.936 & 0.936 & 0.935 & 0.933 & 0.933 \\
 &  & $acl(\hat{c}_{2,2})$ & 0.957 & 0.957 & 0.956 & 0.955 & 0.955 \\\hline
 & 256 & $acl(\hat{c}_{2,1})$ & 0.941 & 0.940 & 0.939 & 0.938 & 0.938 \\
 &  & $acl(\hat{c}_{2,2})$ & 0.955 & 0.955 & 0.954 & 0.954 & 0.953 \\\hline
\end{tabular}
\end{center}
\label{}
\end{table}

\begin{table}[bp]
\caption{$\Sigma=I_p,~\beta=0.01,~eu=0.10$}
\begin{center}
\begin{tabular}{|rr|r|r|r|r|r|r|} \hline
 &  & $p$ & 64 & 128 & 256 & 512 & 1024 \\\hline
$N$ & 64 & $acl(\hat{c}_{2,1})$ & 0.975 & 0.973 & 0.973 & 0.973 & 0.974 \\
 &  & $acl(\hat{c}_{2,2})$ & 0.994 & 0.993 & 0.993 & 0.993 & 0.993 \\\hline
 & 128 & $acl(\hat{c}_{2,1})$ & 0.981 & 0.980 & 0.979 & 0.979 & 0.979 \\
 &  & $acl(\hat{c}_{2,2})$ & 0.993 & 0.992 & 0.992 & 0.993 & 0.992 \\\hline
 & 256 & $acl(\hat{c}_{2,1})$ & 0.984 & 0.984 & 0.983 & 0.983 & 0.982 \\
 &  & $acl(\hat{c}_{2,2})$ & 0.992 & 0.992 & 0.992 & 0.992 & 0.991 \\\hline
\end{tabular}
\end{center}
\label{}
\end{table}

\begin{table}[bp]
\caption{$\rho=0.2,~\beta=0.10,~eu=0.20$}
\begin{center}
\begin{tabular}{|rr|r|r|r|r|r|r|} \hline
 &  & $p$ & 64 & 128 & 256 & 512 & 1024 \\\hline
$N$ & 64 & $acl(\hat{c}_{2,1})$ & 0.887 & 0.885 & 0.885 & 0.886 & 0.883 \\
 &  & $acl(\hat{c}_{2,2})$ & 0.909 & 0.908 & 0.907 & 0.907 & 0.905 \\\hline
 & 128 & $acl(\hat{c}_{2,1})$ & 0.893 & 0.890 & 0.890 & 0.888 & 0.889 \\
 &  & $acl(\hat{c}_{2,2})$ & 0.908 & 0.905 & 0.906 & 0.904 & 0.904 \\\hline
 & 256 & $acl(\hat{c}_{2,1})$ & 0.896 & 0.895 & 0.892 & 0.893 & 0.891 \\
 &  & $acl(\hat{c}_{2,2})$ & 0.906 & 0.905 & 0.903 & 0.904 & 0.902 \\\hline
\end{tabular}
\end{center}
\label{}
\end{table}

\clearpage

\begin{table}[bp]
\caption{$\rho=0.2,~\beta=0.05,~eu=0.20$}
\begin{center}
\begin{tabular}{|rr|r|r|r|r|r|r|} \hline
 &  & $p$ & 64 & 128 & 256 & 512 & 1024 \\\hline
$N$ & 64 & $acl(\hat{c}_{2,1})$ & 0.936 & 0.935 & 0.932 & 0.934 & 0.933 \\
 &  & $acl(\hat{c}_{2,2})$ & 0.958 & 0.957 & 0.955 & 0.956 & 0.955 \\\hline
 & 128 & $acl(\hat{c}_{2,1})$ & 0.940 & 0.936 & 0.940 & 0.939 & 0.940 \\
 &  & $acl(\hat{c}_{2,2})$ & 0.955 & 0.958 & 0.955 & 0.953 & 0.954 \\\hline
 & 256 & $acl(\hat{c}_{2,1})$ & 0.944 & 0.944 & 0.944 & 0.943 & 0.941 \\
 &  & $acl(\hat{c}_{2,2})$ & 0.955 & 0.955 & 0.954 & 0.953 & 0.952 \\\hline
\end{tabular}
\end{center}
\label{}
\end{table}

\begin{table}[bp]
\caption{$\rho=0.2,~\beta=0.01,~eu=0.20$}
\begin{center}
\begin{tabular}{|rr|r|r|r|r|r|r|} \hline
 &  & $p$ & 64 & 128 & 256 & 512 & 1024 \\\hline
$N$ & 64 & $acl(\hat{c}_{2,1})$ & 0.980 & 0.979 & 0.979 & 0.979 & 0.979 \\
 &  & $acl(\hat{c}_{2,2})$ & 0.993 & 0.993 & 0.993 & 0.993 & 0.993 \\\hline
 & 128 & $acl(\hat{c}_{2,1})$ & 0.985 & 0.983 & 0.983 & 0.982 & 0.983 \\
 &  & $acl(\hat{c}_{2,2})$ & 0.993 & 0.992 & 0.992 & 0.992 & 0.992 \\\hline
 & 256 & $acl(\hat{c}_{2,1})$ & 0.986 & 0.985 & 0.985 & 0.985 & 0.985 \\
 &  & $acl(\hat{c}_{2,2})$ & 0.991 & 0.991 & 0.992 & 0.992 & 0.991 \\\hline
\end{tabular}
\end{center}
\label{}
\end{table}

\begin{table}[bp]
\caption{$\rho=0.2,~\beta=0.10,~eu=0.10$}
\begin{center}
\begin{tabular}{|rr|r|r|r|r|r|r|} \hline
 &  & $p$ & 64 & 128 & 256 & 512 & 1024 \\\hline
$N$ & 64 & $acl(\hat{c}_{2,1})$ & 0.885 & 0.881 & 0.879 & 0.881 & 0.877 \\
 &  & $acl(\hat{c}_{2,2})$ & 0.916 & 0.912 & 0.910 & 0.910 & 0.906 \\\hline
 & 128 & $acl(\hat{c}_{2,1})$ & 0.890 & 0.890 & 0.886 & 0.885 & 0.884 \\
 &  & $acl(\hat{c}_{2,2})$ & 0.912 & 0.911 & 0.908 & 0.907 & 0.906 \\\hline
 & 256 & $acl(\hat{c}_{2,1})$ & 0.897 & 0.893 & 0.894 & 0.890 & 0.888 \\
 &  & $acl(\hat{c}_{2,2})$ & 0.912 & 0.909 & 0.909 & 0.906 & 0.904 \\\hline
\end{tabular}
\end{center}
\label{}
\end{table}

\begin{table}[bp]
\caption{$\rho=0.2,~\beta=0.05,~eu=0.10$}
\begin{center}
\begin{tabular}{|rr|r|r|r|r|r|r|} \hline
 &  & $p$ & 64 & 128 & 256 & 512 & 1024 \\\hline
$N$ & 64 & $acl(\hat{c}_{2,1})$ & 0.931 & 0.929 & 0.926 & 0.926 & 0.926 \\
 &  & $acl(\hat{c}_{2,2})$ & 0.962 & 0.961 & 0.958 & 0.957 & 0.956 \\\hline
 & 128 & $acl(\hat{c}_{2,1})$ & 0.939 & 0.936 & 0.934 & 0.933 & 0.932 \\
 &  & $acl(\hat{c}_{2,2})$ & 0.961 & 0.958 & 0.956 & 0.955 & 0.954 \\\hline
 & 256 & $acl(\hat{c}_{2,1})$ & 0.943 & 0.941 & 0.941 & 0.939 & 0.939 \\
 &  & $acl(\hat{c}_{2,2})$ & 0.958 & 0.956 & 0.957 & 0.954 & 0.955 \\\hline
\end{tabular}
\end{center}
\label{}
\end{table}

\begin{table}[bp]
\caption{$\rho=0.2,~\beta=0.01,~eu=0.10$}
\begin{center}
\begin{tabular}{|rr|r|r|r|r|r|r|} \hline
 &  & $p$ & 64 & 128 & 256 & 512 & 1024 \\\hline
$N$ & 64 & $acl(\hat{c}_{2,1})$ & 0.975 & 0.974 & 0.972 & 0.973 & 0.973 \\
 &  & $acl(\hat{c}_{2,2})$ & 0.995 & 0.994 & 0.993 & 0.993 & 0.993 \\\hline
 & 128 & $acl(\hat{c}_{2,1})$ & 0.980 & 0.979 & 0.980 & 0.979 & 0.978 \\
 &  & $acl(\hat{c}_{2,2})$ & 0.994 & 0.993 & 0.993 & 0.992 & 0.992 \\\hline
 & 256 & $acl(\hat{c}_{2,1})$ & 0.984 & 0.984 & 0.984 & 0.984 & 0.982 \\
 &  & $acl(\hat{c}_{2,2})$ & 0.993 & 0.993 & 0.992 & 0.992 & 0.991 \\\hline
\end{tabular}
\end{center}
\label{}
\end{table}

\begin{table}[bp]
\caption{$\rho=0.5,~\beta=0.10,~eu=0.20$}
\begin{center}
\begin{tabular}{|rr|r|r|r|r|r|r|} \hline
 &  & $p$ & 64 & 128 & 256 & 512 & 1024 \\\hline
$N$ & 64 & $acl(\hat{c}_{2,1})$ & 0.902 & 0.893 & 0.890 & 0.888 & 0.887 \\
 &  & $acl(\hat{c}_{2,2})$ & 0.924 & 0.916 & 0.913 & 0.909 & 0.908 \\\hline
 & 128 & $acl(\hat{c}_{2,1})$ & 0.906 & 0.901 & 0.894 & 0.893 & 0.890 \\
 &  & $acl(\hat{c}_{2,2})$ & 0.920 & 0.917 & 0.911 & 0.908 & 0.905 \\\hline
 & 256 & $acl(\hat{c}_{2,1})$ & 0.904 & 0.902 & 0.901 & 0.897 & 0.895 \\
 &  & $acl(\hat{c}_{2,2})$ & 0.915 & 0.913 & 0.912 & 0.908 & 0.906 \\\hline
\end{tabular}
\end{center}
\label{}
\end{table}

\begin{table}[bp]
\caption{$\rho=0.5,~\beta=0.05,~eu=0.20$}
\begin{center}
\begin{tabular}{|rr|r|r|r|r|r|r|} \hline
 &  & $p$ & 64 & 128 & 256 & 512 & 1024 \\\hline
$N$ & 64 & $acl(\hat{c}_{2,1})$ & 0.932 & 0.939 & 0.938 & 0.935 & 0.933 \\
 &  & $acl(\hat{c}_{2,2})$ & 0.954 & 0.962 & 0.960 & 0.958 & 0.956 \\\hline
 & 128 & $acl(\hat{c}_{2,1})$ & 0.938 & 0.945 & 0.941 & 0.943 & 0.939 \\
 &  & $acl(\hat{c}_{2,2})$ & 0.953 & 0.960 & 0.957 & 0.958 & 0.954 \\\hline
 & 256 & $acl(\hat{c}_{2,1})$ & 0.943 & 0.948 & 0.947 & 0.946 & 0.943 \\
 &  & $acl(\hat{c}_{2,2})$ & 0.953 & 0.959 & 0.958 & 0.957 & 0.954 \\\hline
\end{tabular}
\end{center}
\label{}
\end{table}

\begin{table}[bp]
\caption{$\rho=0.5,~\beta=0.01,~eu=0.20$}
\begin{center}
\begin{tabular}{|rr|r|r|r|r|r|r|} \hline
 &  & $p$ & 64 & 128 & 256 & 512 & 1024 \\\hline
$N$ & 64 & $acl(\hat{c}_{2,1})$ & 0.983 & 0.980 & 0.980 & 0.980 & 0.980 \\
 &  & $acl(\hat{c}_{2,2})$ & 0.996 & 0.994 & 0.994 & 0.993 & 0.993 \\\hline
 & 128 & $acl(\hat{c}_{2,1})$ & 0.987 & 0.985 & 0.984 & 0.983 & 0.983 \\
 &  & $acl(\hat{c}_{2,2})$ & 0.995 & 0.994 & 0.993 & 0.992 & 0.992 \\\hline
 & 256 & $acl(\hat{c}_{2,1})$ & 0.988 & 0.987 & 0.987 & 0.985 & 0.986 \\
 &  & $acl(\hat{c}_{2,2})$ & 0.994 & 0.993 & 0.993 & 0.992 & 0.992 \\\hline
\end{tabular}
\end{center}
\label{}
\end{table}

\begin{table}[bp]
\caption{$\rho=0.5,~\beta=0.10,~eu=0.10$}
\begin{center}
\begin{tabular}{|rr|r|r|r|r|r|r|} \hline
 &  & $p$ & 64 & 128 & 256 & 512 & 1024 \\\hline
$N$ & 64 & $acl(\hat{c}_{2,1})$ & 0.897 & 0.889 & 0.881 & 0.880 & 0.881 \\
 &  & $acl(\hat{c}_{2,2})$ & 0.934 & 0.924 & 0.915 & 0.912 & 0.911 \\\hline
 & 128 & $acl(\hat{c}_{2,1})$ & 0.904 & 0.898 & 0.891 & 0.890 & 0.887 \\
 &  & $acl(\hat{c}_{2,2})$ & 0.931 & 0.924 & 0.917 & 0.913 & 0.909 \\\hline
 & 256 & $acl(\hat{c}_{2,1})$ & 0.903 & 0.902 & 0.899 & 0.894 & 0.891 \\
 &  & $acl(\hat{c}_{2,2})$ & 0.923 & 0.921 & 0.917 & 0.912 & 0.908 \\\hline
\end{tabular}
\end{center}
\label{}
\end{table}

\begin{table}[bp]
\caption{$\rho=0.5,~\beta=0.05,~eu=0.10$}
\begin{center}
\begin{tabular}{|rr|r|r|r|r|r|r|} \hline
 &  & $p$ & 64 & 128 & 256 & 512 & 1024 \\\hline
$N$ & 64 & $acl(\hat{c}_{2,1})$ & 0.937 & 0.934 & 0.929 & 0.927 & 0.927 \\
 &  & $acl(\hat{c}_{2,2})$ & 0.974 & 0.969 & 0.963 & 0.960 & 0.958 \\\hline
 & 128 & $acl(\hat{c}_{2,1})$ & 0.945 & 0.942 & 0.937 & 0.936 & 0.932 \\
 &  & $acl(\hat{c}_{2,2})$ & 0.971 & 0.967 & 0.962 & 0.960 & 0.956 \\\hline
 & 256 & $acl(\hat{c}_{2,1})$ & 0.949 & 0.947 & 0.944 & 0.941 & 0.941 \\
 &  & $acl(\hat{c}_{2,2})$ & 0.968 & 0.966 & 0.962 & 0.959 & 0.958 \\\hline
\end{tabular}
\end{center}
\label{}
\end{table}

\begin{table}[bp]
\caption{$\rho=0.5,~\beta=0.01,~eu=0.10$}
\begin{center}
\begin{tabular}{|rr|r|r|r|r|r|r|} \hline
 &  & $p$ & 64 & 128 & 256 & 512 & 1024 \\\hline
$N$ & 64 & $acl(\hat{c}_{2,1})$ & 0.977 & 0.975 & 0.975 & 0.973 & 0.974 \\
 &  & $acl(\hat{c}_{2,2})$ & 0.998 & 0.996 & 0.995 & 0.994 & 0.994 \\\hline
 & 128 & $acl(\hat{c}_{2,1})$ & 0.982 & 0.981 & 0.980 & 0.980 & 0.979 \\
 &  & $acl(\hat{c}_{2,2})$ & 0.997 & 0.996 & 0.995 & 0.994 & 0.993 \\\hline
 & 256 & $acl(\hat{c}_{2,1})$ & 0.986 & 0.985 & 0.985 & 0.984 & 0.983 \\
 &  & $acl(\hat{c}_{2,2})$ & 0.996 & 0.996 & 0.995 & 0.994 & 0.993 \\\hline
\end{tabular}
\end{center}
\label{}
\end{table}

\end{document}